\documentstyle[12pt]{article}
\def\M{{\cal M}}
\def\div{{\rm div}_\mu}
\def\L{{\cal L}_{\xi_X}}
\def\U{{\cal U}(g_c)}
\def\Z{{\cal Z}(g_c)}
\def\A{{\cal A}}
\def\k{k_r}
\def\R{{\bf R}}
\def\C{{\bf C}}

\begin{document}

\bigskip

\begin{center}
{\bf  PSEUDO-K\"AHLER QUANTIZATION ON FLAG MANIFOLDS\\
Alexander V. Karabegov}
\end{center}

\vskip 0.3cm

International Centre for Theoretical Physics, Trieste, Italy
\footnote{On leave of absence from the Joint Institute for Nuclear
Research, LCTA, Dubna 141980, Moscow reg., Russia}

e-mail: karabegov@vxjinr.jinr.ru

\bigskip

\begin{abstract}

A unified approach to geometric, symbol and deformation quantizations on
a generalized flag manifold endowed with an invariant pseudo-K\"ahler
structure is proposed. In particular cases we arrive at Berezin's
quantization via covariant and contravariant symbols.

\end{abstract}

{\bf 1. Introduction}

In the series of papers \cite{C1},\cite{C2},\cite{C3},\cite{C4} a modern
approach to quantization on K\"ahler manifolds was proposed which
combines together geometric quantization \cite{Ko},\cite{S}, symbol
quantization \cite{Ber} and deformation quantization \cite{BF}. The main
idea of this approach can be formulated for quantization on a general
symplectic manifold $\M$ as follows.

\begin{itemize}

\item To give a geometric realization of a family of Hilbert spaces
$H_\hbar$ over the manifold $\M$, parametrized by a small parameter
$\hbar$ which plays a role of Planck constant (by means of geometric 
quantization or of its generalizations).

\item To describe a geometric construction of a symbol mapping from
functions on $\M$ to operators in $H_\hbar$ (the construction of
operator symbols).

\item To choose appropriate algebras ${\cal A}_\hbar$ of symbols
such that the symbol mapping provides a
representation of ${\cal A}_\hbar$ in $H_\hbar$.

\item To find the deformation quantization
which controls the asymptotic expansion of the symbol product as
$\hbar\to 0$ in the same geometric framework.

\end{itemize}

The aim of this paper is to carry out this quantization program on a 
generalized flag manifold, a homogeneous space of a compact semisimple 
Lie group, endowed with an invariant pseudo-K\"ahler structure.

In \cite{K2} it was shown that the theory of spherical Harish-Chandra
modules provides a natural algebraic construction of "mixed" symbol
algebras on a generalized flag manifold, which in particular cases are
algebras of Berezin's covariant and contravariant symbols.

In the present paper we consider an alternative, geometric construction
of such an algebra $\A$ on a flag manifold $\M$, which involves some
invariant pseudo-K\"ahler structure on $\M$. This construction was
independently introduced in \cite{A} and \cite{K3}. Then we describe
a natural geometric representation of the algebra $\A$ in sheaf
cohomology of the quantum line bundle on $\M$.  The description is based
on the Bott--Borel--Weil theorem. It turns out that the algebras of
Berezin's covariant and contravariant symbols correspond to totally
positive and totally negative K\"ahler structures respectively.

Then we consider $\hbar$-parametrized families of symbol algebras
for which the asymptotic expansion of the symbol product as $\hbar\to 0$
leads to deformation quantization with separation of variables
(see \cite{K1}).

{\bf 2. Equivariant families of functions on homogeneous manifolds}

Let $K$ be a real Lie group, $\k$ its real Lie algebra, $\k^*$
the real dual to $\k$. For $X\in\k,\ F\in\k^*$  denote their
pairing by $<F,X>$.

Let $\M$ be a homogeneous $K$-manifold.
Denote by $T_k$ the shift operator by $k\in K$ in
$C^\infty(\M)$, $T_kf(x)=f(k^{-1}x),\ x\in\M,\ f\in C^\infty(\M)$.
We call a family of real smooth functions
$\{f_X\},\ X\in\k,$ on $\M$ a {\it $K$-equivariant family} if
$\k\ni X\mapsto f_X$ is a linear mapping from $\k$ to $C^\infty(\M)$,
$K$-equivariant with respect to the adjoint action
on $\k$ and the shift action on $C^\infty(\M)$, so that
for all $k\in K,\ X\in\k$ holds
$T_k f_X=f_{Ad(k)X}$.

For $X\in\k$ denote by $v_X$ the corresponding fundamental vector field
on $\M$. For $k\in K$ holds the relation
$T_kv_XT_k^{-1}=v_{Ad(k)X},\ X\in\k$, where $v_X$ is treated
as a differential operator in $C^\infty(\M)$.

For a $K$-equivariant family $\{f_X\}$ on
$\M$ and for all $X,Y\in\k$ holds the relation $v_Xf_Y=f_{[X,Y]}$.

Given a $K$-equivariant family $\{f_X\}$ on $\M$,
define a "moment" mapping $\gamma:\M\to\k^*$ $K$-equivariant with
respect to the shift action on $\M$ and coadjoint action on $\k^*$, such
that for all $x\in\M, \ X\in\k$ holds $<\gamma(x),X>=f_X(x)$.

Since $\M$ is homogeneous,  the image of $\gamma$ is
a single coadjoint orbit $\Omega=\gamma(\M)\subset\k^*$.
For $X\in\k$ denote by $v_X^\Omega$ the corresponding
fundamental vector field on $\Omega$. It is Hamiltonian
with respect to the $K$-invariant symplectic structure
given by the Kirillov symplectic form $\omega^\Omega$ on $\Omega$. The
function $f_X^\Omega(F)=<F,X>,\ F\in\Omega,$ is its Hamiltonian function,
i.e., $df_X^\Omega=-i(v_X^\Omega)\omega^\Omega$. For $X,Y\in\k$ holds
$\omega^\Omega(v_X^\Omega,v_Y^\Omega)=f^\Omega_{[X,Y]}$.

{\it Remark.} $K$-equivariant families $\{f_X\}$ on $\M$ are in one-to-one
correspondence with $K$-equivariant mappings $\gamma:\M\to\k^*$.
To a given $\gamma$ there corresponds the coadjoint orbit
$\Omega=\gamma(\M)$ and the family $\{f_X\}$ such that
$f_X=\gamma^*f_X^\Omega=f_X^\Omega\circ\gamma$.

Fix a point $x_0\in\M$ and denote by $K_0\subset K$ the isotropy subgroup
of the point $x_0$. The mapping $\gamma$ and thus the
$K$-equivariant family $\{f_X\}$ itself are completely determined by the
image point $\gamma(x_0)$ which is an arbitrary $K_0$-stable point in
$\k^*$.

For two $K$-equivariant families
$\{f^{(1)}_X\}$ and $\{f^{(2)}_X\}$ their linear combination $\{\alpha
f^{(1)}_X+\beta f^{(2)}_X\}$ is a $K$-equivariant family as well.
Therefore the set of all $K$-equivariant families $\{f_X\}$
on $\M$ is a vector space which can be identified with the subspace
$(\k^*)^{K_0}$ of all $K_0$-stable points in $\k^*$.

Denote by $\omega$ the pullback of the form $\omega^\Omega$ by $\gamma$.
Then $\omega$ is a closed (but not necessarily
nondegenerate) $K$-invariant form on $\M$ such that for $X,Y\in\k$
holds $\omega(v_X,v_Y)=f_{[X,Y]}$ and $df_X=-i(v_X)\omega$.
We say that $\omega$ is associated to the
$K$-equivariant family $\{f_X\}$.

The form $\omega$ is nondegenerate iff the tangent mapping to
$\gamma:\M\to\Omega$ at any point $x\in\M$ is an isomorphism of
the tangent spaces $T_x\M$ and $T_{\gamma(x)}\Omega$ or, equivalently,
if $\gamma$ is a covering mapping.

{\bf 3. Modules of functions on complex homogeneous manifolds}

Let $K$ be a real Lie group with the real Lie algebra $\k$.
Denote by $g_c$ the complexification of $\k$,
$g_c=\k \otimes \C$, by $g_r$ the realification of $g_c$,
and by $J$ the corresponding operator of complex structure in $g_r$,
so that $(g_r,J)$ is isomorphic to $g_c$.

Let the group $K$ act transitively and holomorphically
on a complex manifold $\M$. Then for $X\in \k$ the fundamental
vector field $v_X$ on $\M$ decomposes into the sum of holomorphic
and antiholomorphic vector fields $\xi_X$ and $\eta_X$ respectively,
$v_X=\xi_X+\eta_X$. Therefore $\eta_X=\bar\xi_X$ and for arbitrary
$X,Y\in \k$ $\xi_X$ commutes with $\eta_Y$, $[\xi_X,\xi_Y]=\xi_{[X,Y]}$
and $[\eta_X,\eta_Y]=\eta_{[X,Y]}$. For $X\in\k,\ k\in K$ the following
relations hold, $T_k\xi_XT_k^{-1}=\xi_{Ad(k)X}$ and
$T_k\eta_XT_k^{-1}=\eta_{Ad(k)X}$.

For $Z=X+JY\in g_r,\ X,Y\in\k$,
set $\xi_Z=\xi_X+i\xi_Y$ and $\eta_Z=\eta_X-i\eta_Y$.
Now $(g_r,J)\ni Z\mapsto\xi_Z$ is a $\C$-linear homomorphism
from $(g_r,J)$ to the Lie algebra of holomorphic  vector  fields
on $\M$, and $\eta_Z=\bar\xi_Z$. We get that $g_r$ acts on $\M$ by  real
vector fields $v_Z=\xi_Z+\eta_Z=(\xi_X+\eta_X)+i(\xi_Y-\eta_Y)$
which respect the holomorphic structure on $\M$.

We call a mapping $\k\ni X\mapsto s_X$ of the Lie algebra $\k$ to
$End(C^\infty(\M))$ $K$-equivariant if it is $K$-equivariant with
respect to the adjoint action of $K$ in $\k$ and the shift action in
$C^\infty(\M)$.  This means that $T_k s_XT_k^{-1}=s_{Ad(k)X},\ X\in
\k,k\in K$.

Now we shall define a special $(g_r,K)$-module structure on
$C^\infty(\M)$.  Let $K$ act in $C^\infty(\M)$  by the shifts $T_k,\ k\in
K$, and $g_r$ act by real differential operators $m_Z=v_Z+\varphi_Z,\
Z\in g_r,$ where $\varphi_Z$ is a real smooth function on $\M$, so that
the actions of $K$ and $g_r$ agree in the usual sense.
This means that the actions of the algebra $\k$ as of a subalgebra of
$g_r$ and as of the Lie algebra of $K$ coincide, $m_X=v_X$ for $X\in\k$,
and $T_km_ZT_k^{-1}= m_{Ad(k)Z}$ for $k\in K,\ Z\in g_r$.  In particular,
for $X\in\k$ holds $\varphi_X=0$ and $T_k\varphi_Z=\varphi_{Ad(k)Z}$.
Then we say that there is given an {\it s-module} on $\M$.

For $X,Y\in\k$ set $Z(X,Y)=1/2(X-iJX+Y+iJY)\in g_r\otimes\C$.
The mapping $\k\times\k\ni(X,Y)\mapsto Z(X,Y)$ is a Lie
algebra homomorphism from $\k\times\k$ to $g_r\otimes\C$
(moreover, it extends by $\C$-linearity to an isomorphism of
the complex Lie algebras $(g_r,J)\times (g_r,-J)$ and $g_r\otimes\C$).

For $X\in\k$ introduce a function $f_X=(-1/2)\varphi_{JX}$ on $\M$.
It is easy to check that the functions $f_X,\ X\in\k,$ form a
$K$-equivariant family.  For $X\in\k$ set $l_X=\xi_X+if_X,\
r_X=\eta_X-if_X$.  Notice that $\eta_X=\bar\xi_X$ and that the mappings
$\k\ni X\mapsto l_X$ and $\k\ni X\mapsto r_X$ are $K$-equivariant.

A straightforward calculation shows that
$m_{Z(X,Y)}=l_X+r_Y$, where the mapping $g_r\ni Z\mapsto m_Z$
is extended to $g_r\otimes\C$ by $\C$-linearity.
Taking into account that $(X,0)$ commutes with $(0,Y)$ in $\k\times \k$,
we get the following lemma.

{\bf Lemma 1. } {\it The mappings $\k\ni X\mapsto l_X$
and $\k\ni X\mapsto r_X$ are commuting $K$-equivariant
complex conjugate representations
of $\k$ in $C^\infty(\M)$. }

Suppose there is given a representation of $\k$ in
$C^\infty(\M)$ of the form $\k\ni X\mapsto l_X=\xi_X+if_X$,
where $\{f_X\}$ is a $K$-equivariant family on $\M$
(which is equivalent to $X\mapsto l_X$ being $K$-equivariant).
Then there exists an $s$-module on $\M$ to which the representation
$\k\ni X\mapsto l_X$ is associated.

{\bf Lemma 2. } {\it Let $\k\ni X\mapsto l_X=\xi_X+if_X$
be a $K$-equivariant representation of $\k$ in $C^\infty(\M)$.
For $Z=X+JY\in g_r,\ X,Y\in\k,$ define the function $\varphi_Z=-2f_Y$
on $\M$.  Then the mapping $g_r\ni Z\mapsto m_Z=v_Z+\varphi_Z$ is a
 representation of $g_r$ in $C^\infty(\M)$.  Together with the shift
action of $K$ in $C^\infty(\M)$ it defines an $s$-module on $\M$.

 Proof.} Consider the complex conjugate representation
$\k\ni X\mapsto r_X=\bar\l_X=\eta_X-if_X$ to the representation $X\mapsto
l_X$ of $\k$.  Then for $X,Y\in\k$ and $Z(X,Y)=1/2(X-iJX+Y+iJY)$ we have
as above $m_{Z(X,Y)}=l_X+r_Y$.  In order to show that the mapping $g_r\ni
Z\mapsto m_Z$ is a representation of $g_r$ in $C^\infty(\M)$, it is
enough to show that the representations $X\mapsto l_X$ and $X\mapsto r_X$
of $\k$ commute or, equivalently, that $\xi_X f_Y + \eta_Y f_X=0$.  We
get from the identity $[l_X,l_Y]=l_{[X,Y]}$ that
$f_{[X,Y]}=\xi_X f_Y-\xi_Y f_X$.
Since $v_Y f_X=f_{[Y,X]}$, we get that
$f_{[X,Y]}= -\xi_Y f_X - \eta_Y f_X$. Equating the two expressions
for $f_{[X,Y]}$ we obtain the desired identity.
The rest of the proof is straightforward.

It follows from Lemma 2 that any $s$-module is completely determined
by some $K$-equivariant family of functions $\{f_X\}$ on $\M$ for which
the mapping $\k\ni X\mapsto l_X=\xi_X+if_X$ is a representation of $\k$,
or, equivalently, such that for $X,Y\in\k$ holds the relation
$f_{[X,Y]}=\xi_X f_Y-\xi_Y f_X$. Since this relation is linear with
respect to the family $\{f_X\}$, the set $S$ of all $s$-modules on $\M$
is naturally identified with a linear subspace of the vector space of
$K$-equivariant families of functions on $\M$.

It turns out that one can
give a simple characterization of those $K$-equivariant families of
functions which give rise to $s$-modules. It is given in terms of the
closed 2-form $\omega$ on $\M$ associated to the $K$-equivariant family.

{\bf Theorem 1. } {\it A $K$-equivariant family $\{f_X\}$ on $\M$
corresponds to an $s$-module on $\M$ iff the 2-form $\omega$
associated to $\{f_X\}$ is of the type $(1,1)$ with respect to the
complex structure on $\M$.

Proof.} We have to prove that the form
$\omega$ is of the type $(1,1)$ iff the relation
$f_{[X,Y]}=\xi_X f_Y-\xi_Y f_X$ holds for all $X,Y\in\k$.
We can rewrite this relation in terms of $\omega$ using that
$df_X=-i(v_X)\omega$ and $\omega(v_X,v_Y)=f_{[X,Y]}$ as follows,
$\omega(v_X,v_Y)=\omega(\xi_X,v_Y)-\omega(\xi_Y,v_X)$.
Since $v_X=\xi_X+\eta_Y$, it is equivalent to the relation
\begin{equation}
\omega(\xi_X,\xi_Y)=\omega(\eta_X,\eta_Y)\mbox{ for all } X,Y\in\k.
\end{equation}
If $\omega$ is of the type $(1,1)$, the both
sides of (1) vanish.  Suppose now that (1) is true.
For $Z=X+iY\in g_c=\k\otimes\C,\ X,Y\in\k,$ set
$v_Z=v_X+iv_Y,\ \xi_Z=\xi_X+i\xi_Y$ and $\eta_Z=\eta_X+i\eta_Y$.
It follows from (1) that
\begin{equation}
\omega(\xi_Z,\xi_{Z'})=\omega(\eta_Z,\eta_{Z'})\mbox{ for all } Z,Z'\in
g_c.
\end{equation}

Fix a point $x\in\M$. Since $\M$ is $K$-homogeneous, the vectors
$v_X,\ X\in\k,$ at the point $x$ span the real tangent space
$T_x\M$ and thus the vectors $v_Z,\ Z\in g_c,$ span $T_x\M\otimes\C$.
For arbitrary vectors $\xi,\xi'\in T_x\M\otimes\C$ of the
type $(1,0)$ one can find $Z,Z'\in g_c$ such that $v_Z=\xi$
and $v_{Z'}=\xi'$ at the point $x$. Therefore, at the point $x$,
$\xi_Z=\xi,\ \xi_{Z'}=\xi'$ and $\eta_Z=\eta_{Z'}=0$. It follows from
(2) that $\omega(\xi,\xi')=0$ at the point $x\in\M$ for arbitrary vectors
$\xi,\xi'$ of the type (1,0). Since $\omega$ is real, it follows that it
is of the type (1,1), which completes the proof.

We say that an $s$-module is nondegenerate if the corresponding
2-form $\omega$ is nondegenerate. It follows from Theorem 1 that
the 2-form $\omega$ associated to a nondegenerate $s$-module
is a pseudo-K\"ahler form.

The set of nondegenerate $s$-modules is either empty or it is a dense
open conical (i.e. invariant with respect to the
multiplication by non-zero constants, ${\bf s}\mapsto
t\cdot {\bf s},\ {\bf s}\in S,t\in\R\backslash\{0\}$) subset of $S$.

{\it Example.}  Let $\Omega\subset\k^*$ be a coadjoint orbit
of the group $K$, endowed with an invariant pseudo-K\"ahler polarization.
This means that there is given an invariant complex structure
on $\Omega$ such that the Kirillov form $\omega^\Omega$ is of the type
(1,1). The fundamental vector field $v_X^\Omega,\ X\in\k,$ decomposes
into the sum of a holomorphic and antiholomorphic vector fields
$\xi_X^\Omega$ and $\eta_X^\Omega$ respectively,
$v_X^\Omega=\xi_X^\Omega+\eta_X^\Omega$. Then the functions
$f_X^\Omega,\ X\in\k,$ form a $K$-equivariant family that corresponds to
an $s$-module on $\Omega$.  In particular, the mappings $\k\ni X\mapsto
l_X=\xi_X^\Omega+if_X^\Omega$ and $\k\ni X\mapsto
r_X=\eta_X^\Omega-if_X^\Omega$ are two commuting $K$-equivariant
representations of $\k$ in $C^\infty(\Omega)$.

We are going to associate
to each $s$-module on $\M$ an associative algebra $\A$
whose elements are smooth functions on $\M$.

Extend the representations $\k\ni X\mapsto l_X$ and
$\k\ni X\mapsto r_X=\bar l_X$ to $g_c=\k\otimes\C$
by $\C$-linearity. Then, extending them further to
the universal enveloping algebra $\U$ of $g_c$,
one obtains two commuting $K$-equivariant representations
of $\U$ in $C^\infty(\M)$, $u\mapsto l_u$ and $u\mapsto r_u$, $u\in\U$
($K$ acts on $\U$ by the properly extended adjoint action).

Let $u\mapsto\check u$ denote the standard anti-automorphism of
$\U$ which maps $X\in g_c$ to $-X$.

{\bf Lemma 3. } {\it For $u\in\U$ the following relation holds,
$l_u1=r_{\check u}1$.

Proof.} We prove the Lemma for the monomials
$u_n=X_1\dots X_n,\ X_j\in g_c$, using induction over $n$.  One checks
directly that for $X\in\U$ holds
 $l_X1=if_X=r_{\check X}1$.
Assume that $l_{u_{n-1}}1=r_{\check u_{n-1}}1$
holds. Then $l_{u_{n}}1=
l_{u_{n-1}}l_{X_n}1=
l_{u_{n-1}}r_{\check X_n}1=
r_{\check X_n}l_{u_{n-1}}1=r_{\check X_n}
r_{\check u_{n-1}}1=r_{\check u_n}1$.
The Lemma is proved.

For $u\in\U$ denote $\sigma_u=l_u1$ and let
$\A$ denote the image of the mapping
$\sigma:u\mapsto \sigma_u$ from $\U$
to $C^\infty(\M)$. Notice that the mapping $\sigma$
is $K$-equivariant with respect to the adjoint action
on $\U$ and the action by shifts on $C^\infty(\M)$.

{\bf Lemma 4. } {\it The kernel $I$ of the mapping
$\sigma:\U\to C^\infty(\M)$ is a two-sided ideal in $\U$
and thus $\A$ inherits the algebra structure from the
quotient algebra $\U/I$.

Proof.} It follows from the relation $l_u1=0$
that $I$ is a left ideal, while $r_{\check u}1=0$
shows that $I$ is a right ideal since $u\mapsto\check u$
is an anti-homomorphism. The Lemma is proved.

We shall denote the associative product in $\A$ by $*$.
It follows from Lemma 3 that for $u\in\U,\ f\in\A$
holds $l_uf=\sigma_u*f$ and $r_{\check u}f=
f*\sigma_u$.

{\it Remark.} As a subspace of $C^\infty(\M)$ the algebra $\A$ is the
spherical $(g_r,K)$-submodule of the $s$-module it is associated to,
generated by the constant function 1, which is a spherical
($K$-invariant) vector.

Denote by $\Z$ the center of $\U$. The elements of $\Z$
are stable under the adjoint action of $K$.
Since the mapping $\sigma$ is $K$-equivariant,
$\sigma$ maps the central elements of $\U$ to constants in $\A$.
Thus the restriction of the mapping $\sigma$ to $\Z$
defines a central character $\psi:\Z\to\C$ of the algebra
$\U, \psi(z)=\sigma_z,\ z\in\Z$ (here we identify
the constant functions in $\A$ with the corresponding complex constants).

{\bf 4. Holomorphic differential operators on hermitian line bundles}

Let $\pi:L\to \M$ be a holomorphic hermitian
line bundle over $\M$ with hermitian metrics $h$.
Denote by $L^*$ the bundle $L$ with the zero section removed.
It is a $\C^*$-principal bundle.
A local holomorphic trivialization of $L$ is given by a pair
$(U,s)$ where $U$ is an open chart on $\M$ and
$s:U\to L^*$ is a nonvanishing local holomorphic section of $L$.

We are going to define a pushforward of holomorphic
differential operators on $L$ to the base space $\M$.

A holomorphic differential operator $A$ on $L$
is a global geometric object given locally,
for a holomorphic trivialization $(U_\alpha,s_\alpha)$,
by a holomorphic differential operator $A_\alpha$ on $U_\alpha$.
On the intersection of two charts $U_\alpha$ and  $U_\beta$  the
operators $A_\alpha$ and $A_\beta$  must  satisfy  the  relation
$A_\alpha\varphi_{\alpha\beta}=\varphi_{\alpha\beta}A_\beta$
where  $\varphi_{\alpha\beta}$  is  a   holomorphic   transition
function     on     $U_\alpha\cap     U_\beta$     such     that
$\varphi_{\alpha\beta}s_\alpha=s_\beta$.  In  this  relation  we
consider $\varphi_{\alpha\beta}$ as a multiplication operator.

Holomorphic differential operators on $L$ act on the sheaf of local
holomorphic sections of $L$ and form an algebra.

On each chart $(U_\alpha, s_\alpha)$ introduce a real function
$\Phi_\alpha=-\log h\circ s_\alpha$.

{\bf Lemma 5. } {\it On the intersection of two charts $U_\alpha$ and
$U_\beta$ the following equality holds,
$\Phi_\alpha-\Phi_\beta=\log |\varphi_{\alpha\beta}|^2$.

Proof.} We have $\Phi_\beta=-\log h\circ s_\beta=
-\log h\circ (\varphi_{\alpha\beta}s_\alpha)=
-\log(|\varphi_{\alpha\beta}|^2h\circ s_\alpha)=
\Phi_\alpha-\log |\varphi_{\alpha\beta}|^2$.
The Lemma is proved.

Given a global holomorphic differential operator $A$ on $L$,
consider differential operators
$\check A_\alpha=e^{-\Phi_\alpha}A_\alpha\ e^{\Phi_\alpha}$
on each chart $U_\alpha$

{\bf Lemma 6. } {\it The operators $\check A_\alpha$ define a global
differential operator $\check A$ on $\M$.

Proof.} We have to check that on the intersection of two  charts
$U_\alpha$ and $U_\beta$ holds the equality
$\check  A_\alpha=\check   A_\beta$.
It is equivalent to $e^{-\Phi_\alpha}A_\alpha\  e^{\Phi_\alpha}=
e^{-\Phi_\beta}A_\beta\ e^{\Phi_\beta}$
or $A_\alpha e^{\Phi_\alpha-\Phi_\beta}=
e^{\Phi_\alpha-\Phi_\beta}A_\beta$.
Applying Lemma 5 we get an equivalent equality
$A_\alpha\varphi_{\alpha\beta}\bar\varphi_{\alpha\beta}=
\varphi_{\alpha\beta}\bar\varphi_{\alpha\beta}A_\beta$.
The assertion of the Lemma follows now from the fact that
the holomorphic differential operator $A_\beta$
commutes with the multiplication by the antiholomorphic function
$\bar\varphi_{\alpha\beta}$.

We call $\check A$ the pushforward of the holomorphic differential
operator $A$ on the line bundle $L$ to the base space $\M$.
It is clear that the pushforward mapping $A\mapsto\check A$ is an injective
homomorphism of the algebra of holomorphic differential operators on $L$ into
the algebra of differential operators on $\M$.

Let $\nabla$ denote the canonical holomorphic connection of
the hermitian line bundle $(L,h)$. For a local
holomorphic trivialization $(U_\alpha,s_\alpha)$
a local expression of $\nabla$ on $U_\alpha$ is
$\nabla=d-\partial\Phi_\alpha$.  The curvature $\omega$ of $\nabla$
has a local expression $\omega=i\partial\bar\partial\Phi_\alpha$.

Let the Lie group $K$ act on the line  bundle $\pi:L\to \M$ by
holomorphic line bundle automorphisms which
respect the hermitian metrics $h$.

The metrics $h$ can be considered as a function on $L$,
$L\ni q\mapsto h(q)$.
To each local section $s$ of $L$ over an open set
$U\subset\M$ relate a function $\psi_s$ on
$\pi^{-1}(U)\cap  L^*$ such that $\psi_s(q)q=s(\pi(q))$.
For $t\in\C^*$  holds $h(tq)=|t|^2h(q)$  and
$\psi_s(tq)=t^{-1}\psi_s(q)$.

Any element $X$ of the Lie algebra $k_r$ of $K$ acts  on  $L^*$
by a real vector field $v_X^L$ which is the sum of holomorphic
and antiholomorphic vector fields $\xi_X^L$ and  $\eta_X^L=\bar\xi_X^L$
respectively, $v_X^L=\xi_X^L+\eta_X^L$. The vector fields
$v_X^L,\ \xi_X^L$ and $\eta_X^L$
are homogeneous of order $0$ with respect to the  action
of $\C^*$ on $L^*$. Let $v_X,\ \xi_X$ and $\eta_X$ denote their
projections to $\M$, so $v_X=\xi_X+\eta_X$.

The action of $\xi_X^L$ on the functions $\psi_s$ on $L^*$
can be transferred to the action on the corresponding
local holomorphic sections $s$ of $L$, which defines a
global holomorphic differential operator $A_X$ on $L$.
The object of interest to us will be its
pushforward $\check A_X$ to the base space $\M$.

First, consider a local trivialization of $L^*$ by a local section
$s_0: U\to L^*$, which identifies
$(x,v)\in U\times \C^*$ with $s_0(x)v\in L^*|_U$.
Then, locally, $\xi_X^L=\xi_X-a_Xv\partial/\partial v$ for some
holomorphic function $a_X$ on $U$.
To push forward a holomorphic differential
operator from $L|_U$ to $U$ we use the function
$\Phi=-\log h\circ s_0$ on $U$.
The metrics $h$ at the point $(x,v)\in U\times \C^*$
can be expressed as follows,
$h(x,v)=e^{-\Phi}|v|^2$. Since the metrics $h$ is $K$-invariant, we have
$v_X^Lh=0$. A simple calculation shows then that
\begin{equation}
(\xi_X+\bar \xi_X)\Phi=-a_X-\bar a_X.
\end{equation}
Introduce a function $f_X=-i(a_X+\xi_X\Phi)$ on $U$.
Then (3) means that $f_X$ is real.

{\bf Proposition 1. } {\it The holomorphic differential
operator $A_X$ on $L$ and its pushforward $\check A_X$
to the base space $\M$ can be expressed as follows,
$A_X=\nabla_{\xi_X}+if_X$ and $\check A_X=\xi_X+if_X$.
The mapping $\k\ni X\mapsto\check A_X$ is a $K$
equivariant representation of $\k$ in $C^\infty(\M)$.
The function $f_X$ is globally defined on $\M$
and satisfies the relations
$h(\xi_X^Lh^{-1})=if_X\circ\pi$ and
$df_X=-i(v_X)\omega$.

Proof.} Fix a trivialization of $L$ over an open subset $U\subset\M$,
$L|_U\approx U\times\C$,
and consider a local section of $L$ over $U$, $s:x\mapsto
s(x)=(x,\varphi_s(x))\in U\times \C$.  The function $\psi_s$
corresponding to $s$ is defined by the equality $\psi_s(q)q=s(x)$ for
$q\in L^*,\ x=\pi(q)$.  At the point $q=(x,v)\in U\times\C^*$ we have
$(x,\psi_s(x,v)v)=(x,\varphi_s(x))$, whence
$\psi_s(x,v)=\varphi_s(x)v^{-1}$.  To find the local expression of the
operator $A_X$ calculate its action on $\psi_s(x,v)$,
$(\xi_X-a_Xv\partial/\partial v)(\varphi_sv^{-1})=
(\xi_X\varphi_s+a_X\varphi_s)v^{-1}$. Thus, locally,
$A_X=\xi_X+a_X=(\xi_X-\xi_X\Phi)+(\xi_X\Phi+a_X)=
\nabla_{\xi_X}+if_X$.
Pushing it forward to $U$
we get $\check A_X=e^{-\Phi}(\xi_X+a_X)e^\Phi=
\xi_X+(\xi_X\Phi+a_X)=\xi_X+if_X$.
The $K$-equivariance of the mapping $\k\ni X\mapsto\check A_X$
follows from the fact that $K$ acts on the hermitian line bundle
$(L,h)$ by the line bundle automorphisms which preserve the metrics $h$.
We have locally that $h(\xi^L_Xh^{-1})=
e^{-\Phi}|v|^2((\xi_X-a_Xv\partial/\partial v)e^\Phi|v|^{-2})=
\xi_X\Phi+a_X=if_X\circ\pi$.
To prove the last relation of the Proposition
we notice that $i(\xi_X)\omega$ is of the type $(0,1)$ and
$i(\eta_X)\omega$ is of the type $(1,0)$.We have to show that
$\bar\partial f_X=-i(\xi_X)\omega$ and $\partial f_X=-i(\eta_X)\omega$.
These equalities are complex conjugate, so we prove the former one.
Let  $\{z^k\}$ be local holomorphic coordinates on $U$ and
$\xi_X=a^k(z)\partial/\partial z^k$. Then $\bar\partial f_X=-i\bar\partial
(a_X+\xi_X\Phi)=-i\bar\partial\xi_X\Phi=-ia^k(z)(\partial^2\Phi/\partial z^k
\partial \bar z^l)d\bar z^l$. Taking into account that $\omega=
i(\partial^2\Phi/\partial z^k\partial\bar z^l)dz^k\land d\bar z^l$
we immediately obtain the desired equality, which completes the proof.

It follows from Proposition 1 that
to a hermitian line bundle
$(L,h)\to\M$ on which the group $K$ acts by holomorphic automorphisms
which respect the metrics $h$
there corresponds an $s$-module ${\bf s}$ on $\M$. The relation
$df_X=-i(v_X)\omega$ implies that $\omega(v_X,v_Y)=f_{[X,Y]}$, therefore
the (1,1)-form corresponding to the $s$-module ${\bf s}$ is the curvature
$\omega$ of the canonical connection $\nabla$ on $L$.  The pushforward of
the operator $A_X$ to $\M$ coincides with the operator $l_X$, associated
to ${\bf s},\ \check A_X=l_X$.  The mapping $\k\ni X\mapsto A_X$
can be extended to the homomorphism of the algebra $\U$ to
the algebra of holomorphic differential operators on $L$,
$\U\ni u\mapsto A_u$. Since the pushforward mapping is a
homomorphism of the algebra of holomorphic differential operators on $L$
into the algebra of differential operators on $\M$, we get the following
corollary of Proposition 1.

{\bf Corollary. } {\it To a hermitian line bundle
$(L,h)\to\M$ on which the group $K$ acts by holomorphic automorphisms
which respect the metrics $h$
there corresponds an $s$-module ${\bf s}$ on $\M$ such that
for any $u\in\U$ the pushforward of the
operator $A_u$ from $L$ to $\M$ coincides with the operator $l_u$,
associated to ${\bf s},\ \check A_u=l_u$.}

Denote by $L_{can}$ the canonical line bundle
of $\M$, i.e., the top exterior power of the holomorphic cotangent
bundle ${T^*}'\M$ of $\M$, $L_{can}=\land^m{T^*}'\M,$ where
$m={\rm dim}_\C\M$. Its local
holomorphic sections are the local holomorphic $m$-forms on $\M$.  Let
$\mu$ be a global positive volume form on $\M$.  One can associate to it
a hermitian metrics $h_\mu$ on $L_{can}$ such that for an arbitrary local
holomorphic $m$-form $\alpha$ on $\M$
$h_\mu(\alpha)=\alpha\land\bar\alpha/\mu$.

Recall that the divergence of a vector field $\xi$ with respect to
the volume form $\mu$ is given by the formula
$\div\xi={\cal L}_\xi\mu/\mu$, where ${\cal L}_\xi$ is the Lie derivative
corresponding to $\xi$.

Let the volume form $\mu$ on $\M$ be $K$-invariant.
Then for $X\in\k\ \div v_X=\div\xi_X+\div\eta_X=0$.
Since $\mu$ is real, it follows that $\div\xi_X$ and
$\div\eta_X$ are complex conjugate and thus pure imaginary.
For $X\in\k$ introduce a real function $f_X^{can}=-i\div\xi_X$.

The natural geometric action of $K$ on the hermitian line bundle
$(L_{can},h_\mu)$ by holomorphic line bundle automorphisms
preserves the metrics $h_\mu$.
The corresponding infinitesimal action
of an element $X\in\k$ on the local holomorphic $m$-forms on $\M$ by
the Lie derivative $\L$, defines a global
holomorphic differential operator $A_X$ on  $L_{can}$.

{\bf Proposition 2. } {\it The holomorphic differential operator
$A_X$ on $L_{can}$ is given by the formula $A_X=
\nabla_{\xi_X}+\div\xi_X=\nabla_{\xi_X}+if_X^{can}$.

Proof.} Let $U\subset\M$ be a local coordinate chart with
holomorphic coordinates $\{z^k\}$. The form $\alpha_0=
dz^1\land\dots\land dz^m$ is a local holomorphic trivialization
of $L_{can}$. Set $\Phi=-\log h(\alpha_0)$. Then locally
on $(U,\alpha_0)$ $\mu=e^\Phi\alpha_0\land\bar\alpha_0$ and
$\nabla=d-\partial\Phi$. Since $\xi_X$ is a holomorphic
vector field and $\bar\alpha_0$ is an anti-holomorphic form,
we get $\L\bar\alpha_0=0$. Therefore,
$\L\mu=(\xi_X\Phi)e^\Phi\alpha_0\land\bar\alpha_0
+e^\Phi(\L\alpha_0)\land\bar\alpha_0$. On the other hand,
$\L\mu=(\div\xi_X)\mu=
(\div\xi_X)e^\Phi\alpha_0\land\bar\alpha_0$. Therefore,
$(\div\xi_X)\alpha_0=(\xi_X\Phi)\alpha_0+\L\alpha_0$.
Let $\alpha=f\alpha_0$ be a holomorphic $m$-form on $U$.
The holomorphic function $f$ represents the local holomorphic
section $\alpha$ of $L_{can}$ in the trivialization
$(U,\alpha_0)$. Now $\L\alpha=\L(f\alpha_0)=(\xi_Xf)\alpha_0
+f\L\alpha_0=((\xi_X-\xi_X\Phi)f)\alpha_0+(\div\xi_X)\alpha=
(\nabla_{\xi_X}+\div\xi_X)\alpha$. The Proposition is proved.

Applying Proposition 1 we obtain the following

{\bf Corollary. } {\it The pushforward of the operator $A_X$ to $\M$
is $\check A_X=\xi_X+\div\xi_X=\xi_X+if_X^{can}$.
The mapping $\k\ni X\mapsto \xi_X+if_X^{can}$
is a $K$-equivariant representation of $\k$ in $C^\infty(\M)$.}

This means that if there exists a $K$-invariant measure $\mu$ on $\M$, we
get an $s$-module on $\M$. It is easy to check that if we replace $\mu$
by an arbitrary $K$-invariant measure $c\cdot\mu,\ c\in\R_+,$ we will get
the same functions $f_X^{can},\ X\in\k,$ and thus the same $s$-module. This
$s$-module will be called {\it canonical} and denoted ${\bf s}_{can}$.
If the set of nondegenerate $s$-modules on $\M$ is non-empty, then
there exists $K$-invariant symplectic (pseudo-K\"ahler)
form $\omega$ on $\M$ associated to a nondegenerate $s$-module. 
The corresponding
symplectic volume is $K$-invariant as well, and therefore gives rise to
the canonical $s$-module.

Suppose $\mu$ is a $K$-invariant measure on $\M$, and
$\k\ni X\mapsto l_X=\xi_X+if_X$ is a $K$-equivariant representation of
$\k$ in $C^\infty(\M)$, which corresponds to the $s$-module ${\bf s}\in S$.
For a differential operator $A$ in $C^\infty(\M)$, denote by $A^t$ its
formal transpose with respect to the measure $\mu$, so that for all
$\phi,\psi\in C_0^\infty(\M)$ holds $\int (A\phi)\psi d\mu=
\int \phi(A^t\psi)d\mu$.
Consider the $K$-equivariant representation
$\k\ni X\mapsto (l_{-X})^t=\xi_X-if_X+\div \xi_X$. It corresponds
to the $s$-module which we call dual to ${\bf s}$ and denote by ${\bf s}'$.
Since the canonical module ${\bf s}_{can}$ corresponds to the $K$-equivariant
family $\{-i\div\xi_X\}$, we get ${\bf s}'=-{\bf s}+{\bf s}_{can}$.

{\bf 5. Deformation quantizations with separation of variables}

Recall the definition of deformation quantization on a symplectic manifold
$\M$ introduced in \cite{BF}.

{\it Definition.} Formal differentiable deformation quantization on a
symplectic manifold $\M$ is a structure of associative algebra in the
space of all formal series $C^\infty(\M)[[\nu]]$. The product $\star$
of two elements $f=\sum_{r\geq 0} \nu^r f_r,\ g=\sum_{r\geq 0} \nu^r g_r$
of $C^\infty(\M)[[\nu]]$ is given by the following formula,
\begin{equation}
 f\star g=\sum_r \nu^r \sum_{i+j+k=r} C_i(f_j,g_k),
\end{equation}
where $C_r(\cdot,\cdot),\ r=0,1,\dots$, are bidifferential operators
such that for smooth functions $\varphi,\psi$ on $\M$ holds
$C_0(\varphi,\psi)=\varphi\psi$ and $C_1(\varphi,\psi)-C_1(\psi,\varphi)=
i\{\varphi,\psi\}$. Here $\{\cdot,\cdot\}$ is the Poisson bracket on $\M$,
corresponding to the symplectic structure.

Then the product $\star$ is called a star-product. The star-product
can be extended by the same formula (4) to the space
${\cal F}=C^\infty(\M)[\nu^{-1},\nu]]$
of formal Laurent series with a finite polar part.

Since the star-product is given by bidifferential operators,  it
is localizable, that is, it can be restricted to any open  subset
$U\subset\M$.  For  $U\subset\M$  denote
${\cal F}(U)=C^\infty(U)[\nu^{-1},\nu]]$ and  for
$f,g\in{\cal  F}(U)$
let $L_f$ and $R_g$ denote the left star-multiplication operator
by $f$ and the right  star-multiplication  operator  by  $g$  in
${\cal  F}(U)$  respectively,  so  that  $L_fg=f\star   g=R_gf$.
The operators $L_f$ and $R_g$ commute for all $f,g\in{\cal  F}(U)$.
Let ${\cal L}(U)$ and
${\cal R}(U)$ denote the algebras of left and right star-multiplication
operators in ${\cal F}(U)$  respectively.  It  is  important  to
notice that both left and  right  star-multiplication  operators
are formal Laurent  series  of  differential  operators  with  a
finite polar part (i.e., with finitely many  terms  of  negative
degree of the  formal  parameter  $\nu$).  We call  such
operators formal differential operators.

Let  $\M$ be a complex manifold  endowed with a pseudo-K\"ahler
form $\omega_0$. This means that $\omega_0$ is a real closed
nondegenerate form of the type $(1,1)$. Then $\M$ is a
pseudo-K\"ahler manifold. The form $\omega_0$ defines a symplectic
structure on $\M$.

A formal deformation of pseudo-K\"ahler form $\omega_0$
is a formal series $\omega=\omega_0+\nu\omega_1+\dots$,
where $\omega_r,\ r>0$, are closed, possibly degenerate
forms of the type $(1,1)$ on $\M$. On any contractible
chart $U\subset\M$ there exists a formal potential
$\Phi=\Phi_0+\nu\Phi_1+\dots$ of $\omega$, which means
that $\omega_r=i\partial\bar\partial\Phi_r,\ r\geq 0$.

{\it Definition.} Deformation quantization on a pseudo-K\"ahler manifold
$\M$ is called quantization with separation of variables if for
any open $U\subset\M$ and
any holomorphic function $a(z)$ and antiholomorphic function
$b(\bar z)$ on $U$ left $\star$-multiplication by $a$ and
right $\star$-multiplication by $b$ are point-wise multiplications,
i.e., $L_a=a$ and $R_b=b$.

We call the corresponding $\star$-product a $\star$-product with
separation of variables.

In \cite{K1} a complete description of all deformation quantizations
with separation of variables on an arbitrary  K\"ahler manifold was given.
It was shown that such quantizations are parametrized by the formal
deformations of the original K\"ahler form. The results obtained in
\cite{K1} are trivially valid for pseudo-K\"ahler manifolds as well.

{\bf Theorem 2. } \cite{K1} {\it Deformation quantizations
with separation of
variables on a pseudo-K\"ahler manifold $\M$ are in one-to-one
correspondence with formal deformations of the pseudo-K\"ahler
form $\omega_0$. If there is given a quantization with separation
of variables on $\M$ corresponding to a formal deformation $\omega$
of the form $\omega_0$, $U$ is a contractible coordinate chart on $\M$
with holomorphic coordinates $\{z^k\}$, and $\Phi$ is a formal
potential of $\omega$, then the algebra ${\cal L}(U)$ of the left
$\star$-multiplication operators consists of those formal differential
operators on $U$ which commute with all $\bar z^l$ and
$\partial\Phi/\partial\bar z^l+\nu\partial/\partial\bar z^l$.
Similarly, the algebra ${\cal R}(U)$ of the right $\star$-multiplication
operators on $U$ consists of those formal differential operators
which commute with all $z^k$ and $\partial\Phi/\partial z^k+
\nu\partial/\partial z^k$.}

{\it Remark.} Given the algebra ${\cal L}(U)$, one can recover the
$\star$-product $f\star g$ for $f,g\in{\cal F}(U)$ as follows.  One finds
a unique operator $A\in{\cal L}(U)$ such that $A1=f$. Obviously, $A=L_f$,
whence $f\star g=L_fg$.

Let $({\cal F},\star)$ denote the deformation quantization with separation of
variables on $\M$ corresponding to a formal deformation
$\omega=\omega_0+\nu\omega_1+\dots$ of a pseudo-K\"ahler form  $\omega_0$.
Then for $f,g\in{\cal F}\ f\star g=\sum_r \nu^r C_r(f,g)$ for
bidifferential operators $C_r(\cdot,\cdot)$.
Later we shall meet the product $\tilde\star$ on ${\cal F}$, opposite to
the $\star$-product $\star$. This means that for $f,g\in{\cal F}\
f\tilde\star g=g\star f=\sum_r \nu^r C_r(g,f)$, whence it is
straightforward that $\tilde\star$ is the $\star$-product corresponding
to a formal deformation quantization on the symplectic manifold
$(\M,-\omega_0)$.

Denote by $\tilde{\cal L},\tilde{\cal R}$
the algebras of left and right star-multiplication
operators of the deformation quantization $({\cal F},\tilde\star)$,
and by $\tilde L_f,\tilde R_f$ the operators of
left and right star-multiplication by an element $f\in{\cal F}$
respectively. It is clear that $\tilde L_f=R_f,\tilde R_f=L_f,
\tilde{\cal L}={\cal R},\tilde{\cal R}={\cal L}$. If $a,b$ are,
respectively, a holomorphic and antiholomorphic functions on an open
subset $U\subset\M$, then $\tilde L_b=b$ and $\tilde R_a=a$.
This means that the product $\tilde\star$ is a $\star$-product with
separation of variables on the complex manifold $\bar\M$, opposite
to $\M$ (i.e., with the opposite complex structure).

Let $U$ be a contractible coordinate chart on $\M$
with holomorphic coordinates $\{z^k\}$, and $\Phi$ a formal
potential of $\omega$, then the algebra $\tilde{\cal L}(U)={\cal R}(U)$
consists of formal operators, commuting with all $z^k$ and
$\partial\Phi/\partial z^k+\nu\partial/\partial z^k$. Since on $\bar\M$
holomorphic and antiholomorphic coordinates are swapped, the
formal (1,1)-form on $\bar\M$, corresponding to the quantization
$({\cal F},\tilde\star)$ is $i\bar\partial\partial\Phi=-\omega$.
This (1,1)-form is a formal deformation of the pseudo-K\"ahler form
$-\omega_0$ on $\bar\M$.

Let $\M$ be a $K$-homogeneous complex manifold, ${\bf s}_0,{\bf s}_1,\dots$
be $s$-modules on $\M$, $\{f^n_X\}$ and $\omega_n$
be the $K$-equivariant family and the (1,1)-form on $\M$, respectively,
associated to ${\bf s}_n$.
Then $df^n_X=-i(v_X)\omega_n$.
Assume that ${\bf s}_0$ is nondegenerate, i.e., $\omega_0$ is a
pseudo-K\"ahler form.

Denote by $({\cal F},\star)$ the deformation quantization with
separation of variables on $\M$ corresponding to the formal deformation
$\omega=\omega_0+\nu\omega_1+\dots$ of the pseudo-K\"ahler form  $\omega_0$.
Since all the (1,1)-forms $\omega_n$ are $K$-invariant, the
$\star$-product $\star$ is invariant under $K$-shifts.

For $X\in\k$ denote $f^{(\nu)}_X=f^0_X+\nu f^1_X+\dots$.
Introduce a formal operator
$l^{(\nu)}_X=\xi_X+(i/\nu)f^{(\nu)}_X$.

{\bf Proposition 3. } {\it The mapping $\k\ni X\mapsto l^{(\nu)}_X$
is a Lie algebra homomorphism of $\k$ to the algebra
${\cal L}(\M)$ of the left $\star$-multiplication operators
of the deformation quantization $({\cal F},\star)$. It is
$K$-equivariant with respect to the coadjoint action on $\k$
and the conjugation by shift operators in ${\cal L}(\M)$.

Proof.} The mapping $\k\ni X\mapsto l^{(\nu)}_X$
is a $K$-equivariant Lie algebra homomorphism to the Lie algebra
of formal operators on $\M$ if and only if
$\xi_Xf^{(\nu)}_Y-\xi_Yf^{(\nu)}_X=f^{(\nu)}_{[X,Y]}$ and
$T_kf^{(\nu)}_X=f^{(\nu)}_{Ad(k)X}$ for all
$X,Y\in\k$ and $k\in K$. These relations
follow immediately from the corresponding relations for the functions
$f^n_X$. Theorem 2 tells that in order to show that
$l^{(\nu)}_X\in{\cal L}(\M)$ one has to check that
for a formal potential $\Phi$ of $\omega$  on any contractible
coordinate chart $U$ with holomorphic coordinates $\{z^k\}$
the formal operator $l^{(\nu)}_X=\xi_X+(i/\nu)f^{(\nu)}_X$ commutes with all
$\bar z^l$ and
$\partial\Phi/\partial\bar z^l+\nu\partial/\partial\bar z^l$.
Thus we have to check the equality
\begin{equation}
\xi_X(\partial\Phi/\partial\bar z^l)=i\partial f^{(\nu)}_X/\partial\bar z^l.
\end{equation}
Taking into account that $\omega=i\partial\bar\partial\Phi=
i(\partial^2 \Phi/\partial
z^k\partial\bar z^l)dz^k\land d\bar z^l$ and writing down the local
expression for $\xi_X$, $\xi_X=a^k(z)\partial/\partial z^k$,
we rewrite the left hand side of (5) as follows,
$a^k(z)\partial^2 \Phi/\partial
z^k\partial\bar z^l$. On the other hand,
$i\partial f^{(\nu)}_X/\partial\bar z^l=i<-i(v_X)\omega,
\partial/\partial\bar z^l>=-i\omega(v_X,\partial/\partial\bar z^l)
=a^k(z)\partial^2 \Phi/\partial z^k\partial\bar z^l$, which
proves (5) and completes the proof of the Proposition.

Extend the mapping $\k\ni X\mapsto l^{(\nu)}_X$
to the homomorphism $\U\ni u\mapsto l^{(\nu)}_u$
from $\U$ to ${\cal L}(\M)$ and set $\sigma^{(\nu)}_u
=l^{(\nu)}_u1$.

{\bf Corollary. } {\it The mapping $\U\ni u\mapsto \sigma^{(\nu)}_u$
is a homomorphism from $\U$ to the algebra $({\cal F},\star)$,
$K$-equivariant with respect to the adjoint action on $\U$ and
the shift action on ${\cal F}$.}

It follows that the mapping $\U\ni u\mapsto \sigma^{(\nu)}_u$
maps the elements of the center
$\Z$ of $\U$ to formal series with constant coefficients.

{\bf Lemma 7. } {\it For $z\in \Z$ the operator $l^{(\nu)}_z$ is
scalar and is equal to $\sigma^{(\nu)}_z$.

Proof.} If $A\in{\cal L}(\M)$ then $A=L_f$ for $f=A1\in{\cal F}$.
Therefore $l^{(\nu)}_z=L_{\sigma^{(\nu)}_z}$. Let $B$ denote the
multiplication operator by the formal series with constant coefficients
$\sigma^{(\nu)}_z$. It commutes with all formal differential operators
and therefore $B\in{\cal L}(\M)$. Since $B1=\sigma^{(\nu)}_z$ we get
that $l^{(\nu)}_z=B$. The Lemma is proved.

It was shown in Section 3 that for a given $s$-module ${\bf s}$ on $\M$
the function $\sigma_z=l_z1,\ z\in\Z,$ is scalar and is equal to the
value $\psi(z)$ of the central character $\psi$ associated to ${\bf s}$.
Yet it does not mean that for $z\in Z$ the
corresponding operator $l_z$ is scalar. We shall use deformation
quantization to prove the following proposition.

{\bf Proposition 4. } {\it  Let ${\bf s}_1$ be an arbitrary $s$-module
on $\M$, $l_u,\ u\in\U,$ and $\psi$ be the associated operators and the
central character of $\U$ respectively. If the set of nondegenerate
$s$-modules on $\M$ is non-empty then for $z\in\Z$ holds
$l_z=\psi(z)\cdot{\bf 1}$.}

(We denote by ${\bf 1}$ the identity operator.)

{\it Proof. }  Choose a nondegenerate $s$-module ${\bf s}_0$.
Denote by $\{f_X^j\}$ and by $\omega_j$ the $K$-equivariant
family and the (1,1)-form associated to ${\bf s}_j,\ j=0,1,$ respectively.
Consider a parameter dependent $s$-module ${\bf s}(t)=t{\bf s}_0+{\bf s}_1$.
The $K$-equivariant family $\{f_X\}$ associated
to ${\bf s}(t)$ is such that $f_X=tf_X^0+f_X^1$.
Thus for $X\in\k$ the operator $l_X(t)$
associated to ${\bf s}(t)$ is given by the
formula $l_X(t)=\xi_X+i(tf_X^0+f_X^1)$.
When $t=0$ the operator $l_X(t)$ reduces to the operator
$l_X=\xi_X+if_X^1$ associated to the (possibly degenerate)
$s$-module ${\bf s}_1$.
If we replace the parameter $t$ in $l_X(t)$ by $1/\nu$ we will get
the operator $l^{(\nu)}_X=\xi_X+(i/\nu)(f_X^0+\nu f_X^1)$ of the
deformation quantization with separation of variables $(\cal F,\star)$
which corresponds to the formal (1,1)-form $\omega=\omega_0+\nu\omega_1$.
For $z\in\Z$ the operator $l_z(t)$ is polynomial in $t$. If we
replace $t$ by $1/\nu$ in $l_z(t)$ we will get the operator $l^{(\nu)}_z$
which is scalar by Lemma 7. Therefore $l_z(t)$ is scalar as well. Taking
$t=0$ we get that the operator $l_z(0)=l_z$ associated to the $s$-module
${\bf s}_1$ is scalar. Since $l_z1=\sigma_z=\psi(z)$ it follows that
$l_z=\psi(z)\cdot{\bf 1}$.  This completes the proof.

{\bf Theorem 3. } {\it Let $(L,h)\to\M$ be a hermitian line bundle on $\M$
on which the group $K$ acts by holomorphic automorphisms
which preserve the metrics $h$. The algebra $\U$ acts on $(L,h)$ by
holomorphic differential operators $A_u,\ u\in\U$.
Let ${\bf s}$ be the corresponding $s$-module on $\M$
and $\psi$ be the central character of $\U$ associated to ${\bf s}$.
If the set of nondegenerate $s$-modules on $\M$
is non-empty, the center $\Z$ of $\U$ acts on the sheaf of local
holomorphic sections of $L$ by scalar operators
$A_z=\psi(z)\cdot{\bf 1},\ z\in\Z$.

Proof. } It follows from Proposition 3 and Corollary to Proposition 1
that for $z\in\Z$ the pushforward of the holomorphic differential
operator $A_z$ from  $L$ to $\M$ is scalar and is equal to
$\psi(z)=\sigma_z$.  Now the theorem is a consequence of the fact that
the pushforward mapping $A\mapsto\check A$ is injective.

{\bf 6. $s$-modules on flag manifolds}

We are going to apply the results obtained above to the case of $K$ being
a compact semisimple Lie group. The general facts from the theory of
semisimple Lie groups mentioned below may be found in \cite{W}.

Let $g_c$ be a complex semisimple Lie algebra, $h_c$ its Cartan
subalgebra, $h_c^*$ the dual of $h_c$, $W$ the Weyl group of the pair
$(g_c,h_c)$, $\triangle, \triangle^+, \triangle^-, \Sigma\subset
h_c^*$ the sets of all nonzero, positive, negative and simple roots
respectively, $\delta$ the half-sum of positive roots. For each
$\alpha\in\triangle$ choose weight elements $X_\alpha\in g_c$ such that
$[H_\alpha,X_{\pm\alpha}]= \pm 2X_{\pm\alpha}$ for
$H_\alpha=[X_\alpha,X_{-\alpha}]$.

An element $\lambda\in h_c^*$ is
called dominant if $\lambda(H_\alpha)\geq 0$ for all $\alpha\in\Sigma$,
and is a weight if $\lambda(H_\alpha)\in{\bf Z}$ for all
$\alpha\in\Sigma$.  Denote by ${\cal W}$ the set of all weights
in $h_c^*$ (the weight lattice).

Fix an arbitrary subset $\Theta$ of $\Sigma$ and denote by $<\Theta>$
the set of roots which are linear combinations of elements of $\Theta$.
Then $\Pi=<\Theta>\cup\triangle^-$ is a parabolic subset of $\triangle$.
Denote by $g_c^\Theta$ the Levi subalgebra of $g_c$ generated by $h_c$
and $X_\alpha,\ \alpha\in <\Theta>$, and by $q_c$ the parabolic
subalgebra generated by $h_c$ and $X_\alpha,\ \alpha\in\Pi$.

Denote by $g_r,q_r,g_r^\Theta$ the realifications of $g_c,q_c,g_c^\Theta$
respectively, and by $J$ the complex structure in $g_r$ inherited from
$g_c$.

Let $k_r\subset g_r$ denote the compact form of $g_c$ generated by
$JH_\alpha, X_\alpha-X_{-\alpha},J(X_\alpha+X_{-\alpha}),\
\alpha\in\triangle$.
 Define $k_r^\Theta=k_r\cap g_c^\Theta=k_r\cap q_c$.
It is generated by $JH_\alpha,\ \alpha\in\triangle$, and
$X_\alpha-X_{-\alpha}, J(X_\alpha+X_{-\alpha}),\ \alpha\in<\Theta>$.

Introduce the real Lie algebra $t_r=h_c\cap\k$, the Lie algebra
of a maximal torus in $K$. It is generated by $JH_\alpha,\
\alpha\in\triangle$.

Let $G$ be a complex connected simply connected Lie group with the Lie
algebra $g_r$, $G^\Theta$ and $Q$ the Levi and parabolic subgroups of $G$
with the Lie algebras $g_r^\Theta$ and $q_r$ respectively.
In the rest of this paper $K$ will denote maximal
compact subgroup of $G$ with the Lie algebra $\k$, and $K^\Theta=K\cap
G^\Theta=K\cap Q$. It is known that $K^\Theta$ is the centralizer of a torus
and is connected, and that $G/Q=K/K^\Theta$ is a complex compact
homogeneous manifold (a generalized flag manifold). Denote it by $\M$.

Denote by $x_0$ the class of the unit element of $K$ in $\M$ (the
"origin" of $\M$) and by ${\cal E}$ the set of all
$K^\Theta$-invariant points of $\k$. The set of
$K$-equivariant mappings $\gamma:\M\to\k$ is parametrized by
${\cal E}$ so that $\gamma$ corresponds to $E=\gamma(x_0)\in
{\cal E}$. Since the group $K^\Theta$ is connected, the set
${\cal E}$ is the centralizer of $\k^\Theta$. It is easy to check
that ${\cal E}=\{H\in t_r|\alpha(H)=0$ for all
$\alpha\in<\Theta>\}$.

Denote by $(\cdot,\cdot)$ the Killing form on $g_c$.
It is $\C$-linear, and its restriction to $\k$ is negative-definite.

Identify the dual $\k^*$ of the Lie algebra $\k$ with $\k$
via the Killing form.
We are going to show that any $K$-equivariant mapping
$\gamma:\M\to\k$ (or the $K$-equivariant family defined by $\gamma$)
corresponds to an $s$-module on $\M$.
Let $\Omega\subset\k$ be the orbit of the point
$E=\gamma(x_0)\in{\cal E}$,  $\omega^\Omega$
be the Kirillov 2-form on
$\Omega$, and $v_X^\Omega,\ X\in\k,$ the fundamental vector fields on
$\Omega$. Then the 2-form $\omega$ on $\M$ corresponding to $\gamma$
equals $\gamma^*\omega^\Omega$. It is known that at the point
$E\in\Omega$ for $X,Y\in\k$ holds
$\omega^\Omega(v^\Omega_X,v^\Omega_Y)=(E,[X,Y])$. Thus at the point
$x_0\in\M\ \omega(v_X,v_Y)=(E,[X,Y])$. The tangent space $T_{x_0}\M$
 to the complex manifold $\M$ carries the natural complex structure
 $\tilde J$.  In view of Theorem 1 in order to show that the mapping
 $\gamma$  corresponds
 to an $s$-module it is enough to check that the form $\omega$ on the
 tangent space $T_{x_0}\M$ is of the type (1,1) or, equivalently, that
 for any $v_1,v_2\in T_{x_0}\M$ holds $\omega(v_1,v_2)=\omega(\tilde
 Jv_1,\tilde Jv_2)$.  We can identify $T_{x_0}\M$ as a real vector
 space with the subspace of $\k$ generated by the basis consisting of
the elements $X_\alpha-X_{-\alpha},J(X_\alpha+X_{-\alpha}),\
\alpha\in\triangle^+\backslash<\Theta>$.
Since $g_r/q_r=\k/\k^\Theta$,
 we get that $\tilde J(X_\alpha-X_{-\alpha})=J(X_\alpha+X_{-\alpha})$
 for $\alpha\in\triangle^+\backslash<\Theta>$.
The tangent space
$T_{x_0}\M$ can be represented as the direct sum of 2-dimensional
real subspaces spanned by the vectors
$X_\alpha-X_{-\alpha},J(X_\alpha+X_{-\alpha}),\
\alpha\in\triangle^+\backslash<\Theta>$. These subspaces are mutually
orthogonal with respect to the skew-symmetric form $(E,[\cdot,\cdot])$.
Now, for $\alpha\in\triangle^+\backslash<\Theta>$ we have
$(E,[\tilde J(X_\alpha-X_{-\alpha}),\tilde
JJ(X_\alpha+X_{-\alpha})])=(E,[J(X_\alpha+X_{-\alpha}),
-(X_\alpha-X_{-\alpha})])=
(E,[X_\alpha-X_{-\alpha},J(X_\alpha+X_{-\alpha})])=
i(E,[X_\alpha-X_{-\alpha},X_\alpha+X_{-\alpha}])=2i(E,H_\alpha)$.
Thus the form $\omega$ is of the type (1,1).
For $\alpha\in\triangle\backslash<\Theta>$ the linear functional
${\cal E}\ni H\mapsto\alpha(H)$ is nonzero, therefore the
set ${\cal E}_{reg}=\{H\in{\cal E}|(H,H_\alpha)\neq 0$ for all
$\alpha\in\triangle\backslash<\Theta>\}$ is a dense open subset of
${\cal E}$.  The form $\omega$ is nondegenerate iff $E\in{\cal E}_{reg}$.

It is known that under the adjoint action of the compact group $K$ on
$\k$ the isotropy subgroup of any element of $\k$ is connected.  Now if
$\omega$ is nondegenerate, the isotropy subgroup of $E=\gamma(x_0)$
coincides with $K^\Theta$ and thus the mapping $\gamma:\M\to\Omega$ is a
bijection.

Define a sesquilinear form
$<\cdot,\cdot>$ on $(T\M,\tilde J)$ by the formula
$<v_1,v_2>=\omega(v_1,\tilde Jv_2)-i\omega(v_1,v_2)$. If $\omega$ is
nondegenerate and thus pseudo-K\"ahler, the form $<\cdot,\cdot>$ is the
corresponding pseudo-K\"ahler metrics on $\M$. The vectors
$X_\alpha-X_{-\alpha},\ \alpha\in\triangle^+\backslash<\Theta>$, form a
basis in the complex vector space $(T_{x_0}\M,\tilde J)$. They are
orthogonal with respect to the form $<\cdot,\cdot>$. We have
$<X_\alpha-X_{-\alpha},X_\alpha-X_{-\alpha}>=(E,[X_\alpha-X_{-\alpha},
J(X_\alpha+X_{-\alpha})])=i(E,[X_\alpha-X_{-\alpha},X_\alpha+X_{-\alpha}])=
2i(E,H_\alpha)=2(E,JH_\alpha)$. (Notice that since $E,JH_\alpha\in t_r,\
(E,JH_\alpha)$ is real.) Thus we have proved the following theorem.

{\bf Theorem 4. } {\it To an arbitrary $K$-equivariant mapping
$\gamma:\M\to\k$ there corresponds an $s$-module ${\bf s}$ on $\M$.
It is nondegenerate iff for $E=\gamma(x_0)$ and all
$\alpha\in\triangle\backslash<\Theta>$ holds $(E,H_\alpha)\neq 0$.
The set of nondegenerate $s$-modules on $\M$ is non-empty.
For a nondegenerate ${\bf s}$ the associated mapping
$\gamma:\M\to\Omega=\gamma(\M)$ is a bijection and the pseudo-K\"ahler
structure on $\M$, pushed forward to the orbit $\Omega$ defines a
pseudo-K\"ahler polarization on it.  The index of inertia of the
corresponding pseudo-K\"ahler metrics $<\cdot,\cdot>$ on $\M$ (i.e. the
number of minuses in the signature) equals
$\#\{\alpha\in\triangle^+\backslash<\Theta>|(E,JH_\alpha)<0\}$.}

{\bf 7. Convergent star-products on flag manifolds}

We are going to extend the class of convergent star-products on
generalized flag manifolds introduced in \cite{C2}, using results from
\cite{K2}.  We retain the notations of Section 6. In particular, the
group $K$ is compact semisimple and $\M$ is a generalized flag manifold.

A representation of the group $K$ in a vector space $V$
is called $K$-finite if any vector $v\in V$ is $K$-finite, i.e., the set
$\{kv\},\ k\in K,$ is contained in a finite dimensional subspace of $V$.
If this is the case, $V$ splits into the direct sum of isotypic
components. For a dominant weight $\zeta\in{\cal W}$ denote by
$V^\zeta$ the component isomorphic to a multiple of irreducible
representation of $K$ with highest weight $\zeta$.

For a $K$-homogeneous manifold $M$ denote by $F(M)$ the space of
continuous functions on $M$  $K$-finite with respect to the shift action.
Since $K$ is compact, it follows from the Frobenius theorem that each
isotypic component $F(M)^\zeta$ is finite dimensional.

Let $\Omega\subset\k$ be a $K$-orbit. A function on $\Omega$ is called
regular if it is the restriction of a polynomial function on $\k$. It is
easy to show that the set of all regular functions on
$\Omega$ coincides with $F(\Omega)$ (see, e.g.,\cite{K2}).

Let $d$ be a nonnegative integer. Denote by
${\cal U}_d$ the subspace of $\U$,
generated by all monomials of the form $X_1\dots X_k$, where
$X_1,\dots,X_k\in g_c$ and $k\leq d$. The subspaces
$\{{\cal U}_d\}$ determine the canonical filtration on
$\U$.

The symmetric algebra $S(g_c)$ can be identified with the
space of polynomials on $k_r$, so that the element $X\in g_c$
corresponds to the linear functional on $k_r$,
$\tilde X(Y)=(X,Y),\ Y\in k_r$. Let $S^d(g_c)$
be the space of homogeneous polynomials on $k_r$ of degree $d$.
The graded algebra, associated with the canonical filtration
on $\U$ is canonically isomorphic to $S(g_c)$, so that
${\cal U}_d/{\cal U}_{d-1}$
corresponds to $S^d(g_c)$. For $u\in{\cal U}_d$
let $\underline{u}^{(d)}$ denote the corresponding element of
$S^d(g_c)$.
If $k\leq d$ and $u=X_1\dots X_k\in{\cal U}_d$, then
$\underline{u}^{(d)}=0$ for $k<d$ and
$\underline{u}^{(d)}=\tilde X_1\dots\tilde X_d$ for $k=d$.

We say that a parameter dependent vector $v(\hbar)$ in a vector space $V$
depends rationally on a real parameter $\hbar$ if
$v(\hbar)$ can be represented in a form $v(\hbar)=\sum_j a_j(\hbar) v_j$
for a finite number of elements $v_j\in V$ and rational functions
$a_j(\hbar)$, i.e., $v(\hbar)\in\C((\hbar))\otimes V$, where $\C((\hbar))$
is the field of rational functions of $\hbar$.
Denote by $O(\hbar)\subset\C((\hbar))$
the ring of rational functions of $\hbar$ regular at $\hbar=0$.
Vector $v(\hbar)$ is called regular at $\hbar=0$ if
$v(\hbar)\in O(\hbar)\otimes V$.

Let $v(\hbar)=\sum_r \hbar^r v_r,\ v_r\in V$, be the Laurent expansion of
$v(\hbar)$ at $\hbar=0$. Since $v(\hbar)$  depends rationally on $\hbar$,
its Laurent expansion has a finite polar part. Denote by $\Psi(v(\hbar))$
the corresponding formal Laurent series, $\Psi(v(\hbar))=\sum_r \nu^r
v_r$.

The set $S$ of $s$-modules on $\M$ is a
finite dimensional vector space. Thus we can consider
an $s$-module ${\bf s}(\hbar)$ on $\M$ depending
rationally on $\hbar$ and regular at $\hbar=0$.
Denote by $\omega(\hbar)$ the (1,1)-form associated to
${\bf s}(\hbar)$. It is clear that $\omega(\hbar)$ also depends
rationally on $\hbar$ and is regular at $\hbar=0$. Moreover,
$\Psi({\bf s}(\hbar))=
\sum_{r\geq 0} \nu^r {\bf s}_r$ for some ${\bf s}_r\in S$
and $\Psi(\omega(\hbar))=\sum_{r\geq 0} \nu^r \omega_r$
where $\omega_r$ is the (1,1)-form associated to ${\bf s}_r$.

Denote by $\gamma(\hbar)$, $\gamma_r:\M\to \k$  the $K$-equivariant
mappings and by $\{f_X^{(\hbar)}\}$, $\{f^r_X\}$ the $K$-equivariant
families corresponding to ${\bf s}(\hbar)$, ${\bf s}_r$ respectively.
For $X\in\k$ the function $f^{(\hbar)}_X$ on $\M$ depends rationally on
$\hbar$, is regular at $\hbar=0$ and $\Psi(f^{(\hbar)}_X)=\sum_{r\geq
0}\nu^r f^r_X$.

The $K$-equivariant family $\{(1/\hbar)f^{(\hbar)}_X\}$ corresponds
to the $s$-module\\
$(1/\hbar){\bf s}(\hbar)$.

It will be convenient to us to denote by
$l^{(\hbar)}_u,\ u\in\U,$ the operators on $\M$ associated to
the $s$-module
$(1/\hbar){\bf s}(\hbar)$
(rather than to ${\bf s}(\hbar)$) and set
$\sigma^{(\hbar)}_u=l^{(\hbar)}_u1$.  In particular, for $X\in\k\
l^{(\hbar)}_X=\xi_X+(i/\hbar)f^{(\hbar)}_X$ and
$\sigma^{(\hbar)}_X=(i/\hbar)f^{(\hbar)}_X$.

Let ${\cal D}_\hbar$ denote the algebra of differential operators on $\M$
depending rationally on $\hbar$.

{\bf Lemma 8.} {\it For $u\in{\cal U}_d$
the differential operator $\hbar^d l^{(\hbar)}_u$
belongs to ${\cal D}_\hbar$. It is regular at $\hbar=0$ and
$\lim_{\hbar\to 0}\hbar^dl^{(\hbar)}_u$ is a multiplication operator
by the function $i^d\ \underline{u}^{(d)}\circ\gamma_0$.
In particular, the function
$\hbar^d\sigma^{(\hbar)}_u$ depends rationally on $\hbar$, is
regular at $\hbar=0$
and $\lim_{\hbar\to 0}\hbar^d\sigma^{(\hbar)}_u=
i^d\ \underline{u}^{(d)}\circ\gamma_0$.

Proof. }  The function $f^{(\hbar)}_X$ equals
$f^0_X=\underline{X}^{(1)}\circ\gamma_0$ at $\hbar=0$.
Let $u=X_1\dots X_k,\ X_j\in\k,$ for $k\leq d$. Then $u\in{\cal U}_d$.
We have $\hbar^d l^{(\hbar)}_u= \hbar^d
l^{(\hbar)}_{X_1}\dots l^{(\hbar)}_{X_k}=\hbar^{d-k}
(\hbar\xi_{X_1}+if^{(\hbar)}_{X_1})\dots
(\hbar\xi_{X_k}+if^{(\hbar)}_{X_k})$. Thus the limit
$\lim_{\hbar\to 0}\hbar^d l^{(\hbar)}_u$ equals zero if $k<d$ and equals
$i^d f^0_{X_1}\dots f^0_{X_d}=i^d\ \underline{u}^{(d)}\circ\gamma_0$
if $k=d$, whence the Lemma follows immediately.

For the rest of the section denote $\Omega=\gamma_0(\M)$
and assume that the $s$-module ${\bf s}_0$ is nondegenerate, so that
$\omega_0$ is pseudo-K\"ahler. Then Theorem 4 implies that
the $K$-equivariant mapping $\gamma_0:\M\to\Omega$ is a bijection.
Thus $\gamma_0^*:F(\Omega)\to F(\M)$ is an isomorphism of $K$-modules.

{\bf Proposition 5. } {\it If the $s$-module ${\bf s}_0$ is
nondegenerate, then for any $f\in F(\M)$ there exist elements
$u_j\in{\cal U}_{d(j)}$ for some
numbers $d(j)$ and rational functions $a_j(\hbar)$ regular at $\hbar=0$,
such that $f=\sum_j\hbar^{d(j)}a_j(\hbar)\sigma^{(\hbar)}_{u_j}$ for all
but a finite number of values of $\hbar$.

Proof. }   Fix a dominant weight $\zeta\in{\cal W}$.
The subspace ${\cal U}_d\subset\U$ is invariant under the adjoint action
of the group $K$, and is finite dimensional.
The mapping ${\cal U}_d\ni u\mapsto\underline{u}^{(d)}\in S^d(g_c)$ is
$K$-equivariant, therefore it maps ${\cal U}_d^\zeta$ to
$F(\Omega)^\zeta$.
Since the space $F(\Omega)$ coincides with the space of
regular functions on $\Omega$ and is isomorphic to $F(\M)$,
one can choose elements $u_j\in{\cal U}_{d(j)}^\zeta$ for some
numbers $d(j)$ such that the functions
$f_j=i^{d(j)}\underline{u_j}^{(d(j))}\circ\gamma_0$ form a basis
$\{f_j\}$ in $F(\M)^\zeta$. Since the function
$\tilde f_j=\hbar^{d(j)}\ \sigma^{(\hbar)}_{u_j}$ on $\M$ depends
rationally on $\hbar$ and is regular at $\hbar=0$, the elements of
the matrix $(b_{jk}(\hbar))$ such that
$\tilde f_j=\sum_k b_{jk}(\hbar)f_k$
are rational functions of $\hbar$ regular at $\hbar=0$.  It
follows from Lemma 8 that the matrix  $(b_{jk}(\hbar))$
coincides with the identity matrix at $\hbar=0$. Thus the elements
of the inverse matrix $(a_{kj}(\hbar))=(b_{jk}(\hbar))^{-1}$
such that $f_k=\sum_j a_{kj}(\hbar)\tilde f_j
=\sum_j a_{kj}(\hbar)\hbar^{d(j)}\
\sigma^{(\hbar)}_{u_j}$ are also rational functions of
$\hbar$ regular at $\hbar=0$. Now the Proposition follows from the fact
that the space $F(\M)$ is a direct sum of the subspaces $F(\M)^\zeta$.

Let $(\A_\hbar,*_\hbar)$ denote the algebra of functions on
$\M$ associated to the $s$-module $(1/\hbar){\bf s}(\hbar)$.
Any function $f\in F(\M)$ can be represented in the form
$f=\sum_j\hbar^{d(j)}a_j(\hbar)\sigma^{(\hbar)}_{u_j}$ for some
$u_j\in{\cal U}_{d(j)}$. Thus $f\in\A_\hbar$ for
all but a finite number of values of $\hbar$.

For a function $g\in\A_\hbar$ holds $f *_\hbar g=
\sum_j\hbar^{d(j)}a_j(\hbar)l^{(\hbar)}_{u_j}g$.
We get from Lemma 8 the following corollary to Proposition 5.

{\bf Corollary. } {\it Any functions $f,g\in F(\M)$ are elements of
the algebra $(\A_\hbar,*_\hbar)$ for all but a finite number of values of
$\hbar$. The product $f*_\hbar g$ as a function on $\M$
depends rationally on $\hbar$ and is regular at $\hbar=0$,
i.e., $f *_\hbar g\in O(\hbar)\otimes F(\M)$.}

{\it Remark. } It is easy to show that extending the multiplication
$*_\hbar$ by $O(\hbar)$-linearity we obtain the associative algebra
$(O(\hbar)\otimes F(\M),*_\hbar)$
over the ring $O(\hbar)$ of rational functions of $\hbar$
regular at $\hbar=0$.

Denote by $\omega$ the formal (1,1)-form
$\Psi(\omega(\hbar))=\omega_0+\nu\omega_1+\dots$. It is a formal
deformation of the pseudo-K\"ahler form $\omega_0$. Denote by $({\cal
F},\star)$ the deformation quantization with separation of variables on
$\M$ corresponding to $\omega$.

Set $f^{(\nu)}_X=\Psi(f^{(\hbar)}_X)=
f^0_X+\nu f^1_X+\dots$. Then $\Psi(l^{(\hbar)}_X)=
\xi_X+(i/\nu)f^{(\nu)}_X$.
It follows from Proposition 3 that for $X\in\k$ the operator
$l^{(\nu)}_X=\xi_X+(i/\nu)f^{(\nu)}_X$ belongs to the algebra
${\cal L}$ of left $\star$-multiplication operators of the deformation
quantization $({\cal F},\star)$. It is easy to check that
the mapping ${\cal D}_\hbar\ni A\mapsto\Psi(A)$
is a homomorphism from ${\cal D}_\hbar$ to the algebra of formal
differential operators on $\M$, therefore for $u\in\U$
$\Psi(l^{(\hbar)}_u)=l^{(\nu)}_u\in{\cal L}$.

Represent a function $f\in F(\M)$ in the form
$f=\sum_j\hbar^{d(j)}a_j(\hbar)\sigma^{(\hbar)}_{u_j}$ for some
$u_j\in{\cal U}_{d(j)}$ and consider the operator $A=
\sum_j\hbar^{d(j)}a_j(\hbar)l^{(\hbar)}_{u_j}\in{\cal D}_\hbar$.
It follows from Lemma 8 that $A$ is regular at $\hbar=0$.
It is straightforward that $\Psi(A)\in{\cal L}$ and $A1=f$, whence
one can easily obtain that $\Psi(A)1=f$ and therefore $\Psi(A)=L_f$.

For $g\in F(\M)$ the product $f*_\hbar g=Ag$ is a
function on $\M$ which depends rationally on $\hbar$ and is regular at
$\hbar=0$.  Therefore the product $f*_\hbar g$
expands to the uniformly and absolutely convergent Taylor series in
$\hbar$ at $\hbar=0$. Finally, $\Psi(f*_\hbar g)=\Psi(Ag)=L_fg=f\star g$.
Thus we have proved the following theorem.

{\bf Theorem 5. } {\it Let ${\bf s}(\hbar)$ be an $s$-module on $\M$ which
depends rationally on the parameter $\hbar$ and is regular
at $\hbar=0$, and $\omega(\hbar)$ be the associated (1,1)-form.
Then $\Psi({\bf s}(\hbar))=
\sum_{r\geq 0} \nu^r {\bf s}_r$ for some ${\bf s}_r\in S$.
Assume that the $s$-module ${\bf s}_0$ is nondegenerate
and denote by $(\A_\hbar,*_\hbar)$ the algebra of functions associated
to the $s$-module $(1/\hbar){\bf s}(\hbar)$.
Any functions $f,g\in F(\M)$ belong to $\A_\hbar$ for all but a finite
number of values of $\hbar$. The product $f*_\hbar g$ expands to the
uniformly and absolutely convergent Taylor series in $\hbar$ at the point
$\hbar=0$,\\
$f*_\hbar g=\sum_{r\geq 0} \hbar^r C_r(f,g)$, where
$C_r(\cdot,\cdot), r=0,1,\dots$, are bidifferential operators which
define the deformation quantization with separation of variables on $\M$
corresponding to the formal deformation
$\omega=\Psi(\omega(\hbar))=\omega_0+\nu\omega_1+\dots$ of the
pseudo-K\"ahler form $\omega_0$.}

{\bf 8. Characters associated to $s$-modules on flag manifolds}

Since the group $K$ is compact, there exists the
$K$-invariant measure $\mu$ of the total volume 1 on the flag
manifold $\M$.  Let ${\bf s}$ be an $s$-module on $\M$ and $\A$ the
corresponding algebra of functions on $\M$.

It is known that  $\U= \Z\oplus[\U,\U]$ (see \cite{W}). Let  $\U\ni
u\mapsto u^0$ denote the corresponding projection of  $\U$ onto  $\Z$.
Recall the following

{\it Definition} \ (see \cite{W}). A linear form  $\kappa:\U\to
{\rm C}$ is called a character of $g_c$ if:

(1) $\kappa(uv)=\kappa(vu),\ \kappa(1)=1;$

(2) $\kappa(u^0v)=\kappa(u)\kappa(v)\quad
(u,v\in {\cal U}(g_c)).$

Thus one has $\kappa(u)=\kappa(u^0)$ for all  $u\in\U$.
Moreover, $\kappa$ is then a homomorphism of  $\Z$ into ${\rm
C}$, a central character of $g_c$.  This central character determines
$\kappa$ completely.

Fix an $s$-module {\bf s} on $\M$ and let $l_u,\ u\in\U,$ be the
operators on $\M$, associated to {\bf s}, $\sigma_u=l_u1$ and $\psi$
be the corresponding central character of
$\U,\ \psi(z)=\sigma_z,\ z\in\Z$.

{\bf Proposition 6. } {\it  A linear form
$\kappa(u)=\int_\M \sigma_u\ d\mu,\ u\in\U$, on $\U$
is a character of $g_c$. For  $z\in\Z\quad \kappa(z)= \psi(z)$.}

{\it Proof.} Since the mapping $u\mapsto
\sigma_u=l_u1$ is $K$-equivariant and the measure $\mu$ is $K$-invariant,
for $k\in K$ holds $\kappa(Ad(k)u)=\kappa(u)$ or,
infinitesimally, for $X\in k_r\ \kappa(Xu-uX)=0$, therefore
$\kappa(uv)=\kappa(vu)$.  The measure $\mu$ is of the total
volume 1 and for  $z\in\Z\quad \sigma_z=\psi(z)$ is scalar, therefore
$\kappa(z)=\psi(z)$. In particular, $\kappa(1)=1$.  Thus (1) is
proved. Now, using that for $u\in\U$ holds $u^0\in\Z$
and $l_{u^0}=\sigma_{u^0}$ is scalar, we get
$\kappa(u^0v)=\int l_{u^0v}1\ d\mu=
\int l_{u^0}l_v1\ d\mu=\sigma_{u^0}
\int l_v1\ d\mu=\kappa(u^0)\kappa(v)$.
This completes the proof of the Proposition.

{\it Remark. } Proposition 6 implies that to any $s$-module ${\bf s}$ on
$\M$ there corresponds a character $\kappa$ of the Lie algebra $g_c$. The
central character $\Z\ni z\mapsto\kappa(z)$ of $\U$ coincides with the
central character $\psi$ associated to ${\bf s}$. Moreover, the mapping
$\A\ni f\mapsto t(f)=\int_\M f\ d\mu$ is a trace on the algebra $\A$,
i.e.,  $t(f*g)=t(g*f)$ for $f,g\in\A$.

Let $\tau$ be an $n$-dimensional irreducible representation of $g_c$.
Since for  $z\in\Z\quad \tau(z)$ is
scalar, it is straightforward that
$\kappa_\tau=(1/n){\rm tr}\ \tau$
is a character of $g_c$. In particular,
$\Z\ni z\mapsto\kappa_\tau(z)$ is the central character
of the representation $\tau$.

{\bf Proposition 7. } {\it Let $\tau$ be an $n$-dimensional irreducible
representation of $\U$ in the vector space $V$ and ${\bf s}$ be an $s$-module
on $\M$.  If the central character $\psi$ associated to
${\bf s}$ coincides with the central character of the representation $\tau$,
then there exists a representation $\rho$ of the algebra $\A$ in the same
vector space $V$ such that $\tau=\rho\circ\sigma$.

Proof. } Since the characters $\kappa_\tau=(1/n){\rm tr}\ \tau$
and $\kappa=\int_\M \sigma\ d\mu$ of $g_c$ coincide on $\Z$,
they coincide identically.
We have that for $u,v\in\U$
\begin{equation}
{\rm tr}(\tau(u)\tau(v))={\rm tr}\ \tau(uv)=n\int \sigma_{uv}\
d\mu=n\int \sigma_u * \sigma_v\ d\mu.
\end{equation}
Assume that $\sigma_v=0$. Then the last expression in (6) is zero for all
$u\in\U$.  Since $\tau$ is irreducible, $\tau(u)$ is an arbitrary
endomorphism of the representation space $V$, therefore  $\tau(v)=0$.
Thus the representation $\tau$ factors through the mapping
$\sigma:\U\to\A,\ \tau=\rho\circ\sigma$.
The Proposition is proved.

{\bf 9. Holomorphic line bundles on flag manifolds}

The Levi subgroup $G^\Theta\subset G$ is reductive.
The pair $(g_c^\Theta,h_c)$ has a root system $<\Theta>$. Induce the
ordering on $<\Theta>$ from $\triangle$.
Denote by $W^\Theta$ the Weyl group of the pair $(g_c^\Theta,h_c)$,
by $\delta_\Theta$ the half-sum of positive roots from $<\Theta>$,
and set $\delta'_\Theta=\delta-\delta_\Theta$.
The one-dimensional holomorphic representations (the holomorphic
characters) of $G^\Theta$ are parametrized by the set ${\cal W}^\Theta$
of $W^\Theta$-invariant weights from ${\cal W}$.  The parabolic group $Q$
is a semi-direct product of $G^\Theta$ and of the unipotent radical $R$
of $Q$. For $\lambda\in{\cal W}^\Theta$ denote by $\chi_\lambda$ the
holomorphic character of $Q$, which is trivial on $R$, and whose
restriction to $G^\Theta$ is the character of $G^\Theta$ parametrized by
$\lambda$. For $H\in h_c\ \chi_\lambda(\exp H)=\exp\lambda(H)$.

Denote by $\C_\lambda$ a one-dimensional complex vector space with the
action of $Q$ given by $\chi_\lambda$. Consider the holomorphic line
bundle $L_\lambda=G\times_Q {\C}_\lambda$. It is the coset space of
$G\times\C_\lambda$ under the equivalence $(gq,v)=(g,\chi_\lambda(q)v),\
g\in G,q\in Q,v\in\C_\lambda$. The group $G$ acts on $L_\lambda$
as follows, $G\ni g_0:(g,v)\mapsto (g_0g,v)$. Since $G/Q=K/K^\Theta$, one
has an alternative description of $L_\lambda,\
L_\lambda=K\times_{K^\Theta}\C_\lambda$. Using that description one can
define a $K$-invariant hermitian metrics $h$ on $L_\lambda$ setting
$h(k,v)=|v|^2$. It follows from Iwasawa decomposition that each element
$g\in G$ can be (non-uniquely) represented as a product $g=kq$ for some
$k\in K,q\in Q$.  Thus for $g=kq$ we get
$h(g,v)=h(kq,v)=h(k,\chi_\lambda(q)v)=|\chi_\lambda(q)v|^2$.

It follows from the results obtained in  Sect. 4, that to the hermitian
line bundle $(L_\lambda,h)$ there corresponds an $s$-module on $\M$.
Denote it ${\bf s}_\lambda$.
Let $\{f^\lambda_X\},\ X\in\k,$
be the corresponding $K$-equivariant family which
defines the mapping $\gamma:\M\to\k$ such that $(\gamma(x),X)=f^\lambda_X(x)$
for all $x\in\M,\ X\in\k$.  We are
going to apply Theorem 4 to the $s$-modules
${\bf s}_\lambda,\ \lambda\in{\cal W}^\Theta$.
Calculate the element $E^\lambda=\gamma(x_0)$.
Since $E^\lambda\in t_r$, in order to determine $E^\lambda$ it is
enough to consider only the pairing  $(E^\lambda,JH_\alpha)$ for all
$\alpha\in\triangle$.
For $Z\in g_r$ let $v_Z^L$ be the fundamental vector
field on $L_\lambda$, then $\xi_Z^L=(1/2)(v_Z^L-iv_{JZ}^L)$ is its
holomorphic component. For $\varphi\in C^\infty(L_\lambda)\
v_Z\varphi(g,v)=(d/dt)\varphi(\exp(-tZ)g,v)|_{t=0},\ g\in
G,v\in\C_\lambda$.

Using Proposition 1 and taking into
account $K$-invariance of the metrics $h$ we get
$if^\lambda_X\circ\pi=h(\xi_X^Lh^{-1})=(-i/2)h(v_{JX}^Lh^{-1})$ for
$X\in\k$.  Thus $f^\lambda_{JH_\alpha}(x_0)=(1/2)h(v_{H_\alpha}h^{-1})=
(1/2)h(e,v)(d/dt)(h(\exp (-tH_\alpha),v))^{-1}|_{t=0}=
(1/2)(d/dt)(|\chi_\lambda(\exp (-tH_\alpha))|^{-2})|_{t=0}=
(1/2)(d/dt)(|\exp \lambda(tH_\alpha)|^2)|_{t=0}=\lambda(H_\alpha)$.
Now $E^\lambda=\gamma(x_0)$ is the element of $t_r$ such that
$(E^\lambda,JH_\alpha)=\lambda(H_\alpha)$ for all $\alpha\in\triangle$.
The following proposition is a direct consequence of Theorem 4.

{\bf Proposition 8. } {\it  The $s$-module ${\bf s}_\lambda$
corresponding to the holomorphic hermitian line bundle $(L_\lambda,h),\
\lambda\in{\cal W}^\Theta,$ is nondegenerate iff for all
$\alpha\in\triangle\backslash<\Theta>$ holds $\lambda(H_\alpha)\neq 0$.
In this case the index of inertia of the corresponding
pseudo-K\"ahler metrics on $\M$ equals
$\#\{\alpha\in\triangle^+\backslash<\Theta>|\lambda(H_\alpha)<0\}$.}

{\bf Lemma 9. } {\it The canonical line bundle $L_{can}$ on $\M$
is isomorphic to the bundle $L_\lambda$ for $\lambda=-2\delta'_\Theta$.
The canonical $s$-module ${\bf s}_{can}$ on $\M$ coincides with
${\bf s}_\lambda$ for $\lambda=-2\delta'_\Theta$.

Proof. } The isotropy subgroup $Q\subset G$ of the
point $x_0\in\M=G/Q$ acts on the fibers of $G$-bundles at $x_0$.
The fiber of the line bundle $L_\lambda$ at $x_0$ is isomorphic as a
$Q$-module to $\C_\lambda$.  On the other hand, the holomorphic tangent
space of $\M$ at $x_0,\ T'_{x_0}\M,$ is isomorphic as a $Q$-module to
$g_c/q_c$ under the adjoint action. For $H\in h_c$ the operator $ad(H)$
on $g_c/q_c$ is diagonal in the basis $\{X_\alpha+q_c\},\
\alpha\in\triangle\backslash\Pi,$ and takes the eigenvalue $\alpha(H)$ on
$X_\alpha+q_c$. The element $H\in h_c$ acts
on $\land^m (g_c/q_c)$ by the scalar $2\delta'_\Theta(H)$, where $m={\rm
dim}_\C\M$.  Therefore the element $q\in Q$ acts on $\land^m (g_c/q_c)$
by $\chi_\lambda(q)$ for $\lambda=2\delta'_\Theta$. The Lemma follows
from the fact that the fiber of the canonical line bundle $L_{can}$ at
$x_0$ is dual to $\land^m (g_c/q_c)$ as a $Q$-module.

Now we shall use a particular case of the Bott-Borel-Weil theorem
concerning cohomological realizations of finite dimensional irreducible
holomorphic representations of the group $G$ in the sheaf
cohomologies of line bundles over $\M=G/Q$ (see \cite{Kn}). Let
$H^i(\M,{\cal S}L_\lambda)$ denote the space of $i$-dimensional
cohomology with coefficients in the sheaf of germs of holomorphic
sections of the line bundle $L_\lambda$. The action of the group $G$ on
$L_\lambda$ gives rise to the action of $G$ on the local holomorphic
sections of $L_\lambda$, which induces the action of $G$ in the
cohomology spaces $H^i(\M,{\cal S}L_\lambda)$.

{\bf Theorem 6. } (Bott-Borel-Weil) {\it  Let
$\lambda\in{\cal W}^\Theta,\
k=\#\{\alpha\in\triangle^+|(\lambda+\delta)(H_\alpha)<0\}$.
 If $(\lambda+\delta)(H_\alpha)=0$
for some $\alpha\in\triangle$ then $H^i(\M,{\cal S}L_\lambda)=0$ for all
$i$.  If $(\lambda+\delta)(H_\alpha)\neq 0$ for all $\alpha\in\triangle$
one can choose $w\in W$ so that $w(\lambda+\delta)$ is dominant.
Then $\zeta=w(\lambda+\delta)-\delta$ is dominant as well. For all
$i\neq k$ $H^i(\M,{\cal S}L_\lambda)=0$. The representation of
the group $G$ in $H^k(\M,{\cal S}L_\lambda)$ is isomorphic to the
irreducible finite dimensional holomorphic representation of $G$ with
highest weight $\zeta$.}

Assume that an irreducible finite dimensional holomorphic representation
$\tau$ of the group $G$ is realized in the cohomology space
$H^k(\M,{\cal S}L_\lambda)$ as in Theorem 6. Retain the same notation
for the representations of the Lie algebra $g_c$ and of its universal
enveloping algebra $\U$ which correspond to $\tau$.
The action of the Lie algebra $g_c$ on $L_\lambda$ by holomorphic
differential operators can be extended to the action of $\U$.
For $u\in\U$ denote by $A_u$ the corresponding holomorphic differential
operator on $L_\lambda$. It induces the representation operator
$\tau(u)$ in $H^k(\M,{\cal S}L_\lambda)$.

According to Theorem 3, for $z\in\Z$ the
holomorphic operator $A_z$ on $L_\lambda$ is scalar and is equal to the
value $\psi(z)$ of the central character $\psi$ associated to the
$s$-module ${\bf s}_\lambda$. It follows
immediately that the central character of the representation
$\tau$ coincides with $\psi$.  As in the proof of
Proposition 7 we obtain that for $u\in\U$ holds the
equality
\begin{equation}
n\int_\M\sigma_u\ d\mu={\rm tr}\ \tau(u),
\end{equation}
where $n={\rm dim}\
\tau$. According to Proposition 7, there exists
a representation $\rho$ of the algebra $\A$ of functions on $\M$
associated to the $s$-module ${\bf s}_\lambda$ in the space
$H^k(\M,{\cal S}L_\lambda)$, such that $\tau=\rho\circ\sigma$.

Let $\{f^\lambda_X\},\ X\in\k,$
be the  $K$-equivariant family associated to ${\bf s}_\lambda$.
Since $\sigma_X=if^\lambda_X$ for $X\in\k$,
the algebra $\A$ contains the functions
$f^\lambda_X,\ X\in\k$, and is generated by them.
The algebra $\k$ acts on $L_\lambda$ by the holomorphic differential
operators $A_X=\nabla_{\xi_X}+if^\lambda_X,\ X\in\k,$ due to Proposition
1. The operator $A_X$ induces in $H^k(\M,{\cal S}L_\lambda)$ the
representation operator $\rho(if^\lambda_X)$.

{\it Remark. } For $X\in\k$ consider the operator
$Q_X=\nabla_{v_X}+if^\lambda_X$. The operator $A_X$ differs from $Q_X$
by the anti-holomorphic operator $\nabla_{\eta_X}$, which annihilates the
local holomorphic sections of $L_\lambda$ and thus induces the trivial
action on the sheaf cohomology. Therefore the operator $Q_X$ also
induces in $H^k(\M,{\cal S}L_\lambda)$ the operator $\rho(if^\lambda_X)$.
If the $s$-module ${\bf s}_\lambda$ is nondegenerate,
the curvature form $\omega$ of the connection $\nabla$ is
symplectic and the function $f^\lambda_X$ is a
Hamiltonian of the fundamental vector field $v_X$ on $\M$.
Then the operator $Q_X$ is the operator of geometric quantization
corresponding to the function $f^\lambda_X$.

We see that the Bott-Borel-Weil theorem provides a natural geometric
representation of the algebra $\A$ in the sheaf cohomology space of the
line bundle $L_\lambda$.

{\bf Theorem 7. } {\it  Let $\lambda\in{\cal W}^\Theta$ be such that
$(\lambda+\delta)(H_\alpha)\neq 0$ for all $\alpha\in\triangle$,
$\A$ and $\{f^\lambda_X\},\ X\in\k,$ be the algebra of functions on
$\M$ and the $K$-equivariant family associated to the
$s$-module ${\bf s}_\lambda$ respectively. The algebra $\A$ is generated by its
elements $f^\lambda_X,\ X\in\k$. Set
$k=\#\{\alpha\in\triangle^+|(\lambda+\delta)(H_\alpha)<0\}$.  There exists
a unique finite dimensional irreducible representation $\rho$ of the
algebra $\A$ in the space $H^k(\M,{\cal S}L_\lambda)$ such that for all
$X\in\k$ the representation operator $\rho(if^\lambda_X)$ is induced from
the holomorphic differential operator $\nabla_{\xi_X}+if^\lambda_X$ on
$L_\lambda$.  There exists an element $w\in W$ such that
$\zeta=w(\lambda+\delta)-\delta$ is a dominant weight of the Lie
algebra $g_c$. The representation $\tau=\rho\circ\sigma$ of $g_c$
in $H^k(\M,{\cal S}L_\lambda)$ is irreducible with highest weight $\zeta$.}

Denote by $w_0$ and $w_0^\Theta$ the elements of the maximal reduced
length in the Weyl groups $W$ and $W^\Theta$ respectively.
Let $\tau$ be the irreducible finite dimensional representation of the
algebra $g_c$ with highest weight $\zeta$. It is known that the dual
representation $\tau'$ has the highest weight $\zeta'=-w_0\zeta$.

{\bf Lemma 10. } {\it Let $\lambda\in{\cal W}^\Theta$ and $w\in W$
be such that $\zeta=w(\lambda+\delta)-\delta$ is a dominant weight.
Then $\lambda'=-\lambda-2\delta'_\Theta\in{\cal W}^\Theta$ and
there exists an element $w'\in W$ such that
$\zeta'=w'(\lambda'+\delta)-\delta$.
If $(\lambda+\delta)(H_\alpha)\neq 0$ for all $\alpha\in\triangle$
and $k=\#\{\alpha\in\triangle^+|(\lambda+\delta)(H_\alpha)<0\}$
then $(\lambda'+\delta)(H_\alpha)\neq 0$
for all $\alpha\in\triangle$ and $\#\{\alpha\in\triangle^+|(\lambda'
+\delta)(H_\alpha)<0\}=m-k$, where $m={\rm dim}_\C\M$.

Proof. } For $\alpha\in\Theta$ the reflection $s_\alpha\in W^\Theta$
maps $\alpha$ to $-\alpha$ and preserves both
$\triangle^+\backslash\{\alpha\}$ and $<\Theta>$. It follows that
the group $W^\Theta$ preserves the set
$\triangle^+\backslash<\Theta>$, whence  $-2\delta'_\Theta\in{\cal
W}^\Theta$ and therefore $\lambda'=-\lambda-2\delta'_\Theta\in{\cal
W}^\Theta$. The element $w_0^\Theta$ maps $<\Theta>^+$ to $<\Theta>^-$
and preserves $\triangle^+\backslash<\Theta>$, whence
$w_0^\Theta\delta_\Theta=-\delta_\Theta$.
Take $w'=w_0ww_0^\Theta$, then
$w'(\lambda'+\delta)=
w'(-\lambda-\delta'_\Theta+\delta_\Theta)=
w_0ww_0^\Theta(-\lambda-\delta'_\Theta+\delta_\Theta)=
w_0w(-\lambda-\delta'_\Theta-\delta_\Theta)=
w_0w(-\lambda-\delta)=w_0(-\zeta-\delta)=
-w_0\zeta+\delta=\zeta'+\delta$, thus
$\zeta'=w'(\lambda'+\delta)-\delta$.
For $\alpha\in\triangle$ we have $(\lambda'+\delta)(H_\alpha)=
(-\lambda-\delta'_\Theta+\delta_\Theta)(H_\alpha)=
(w_0^\Theta(-\lambda-\delta'_\Theta+\delta_\Theta))
(H_{w_0^\Theta\alpha})=
(-\lambda-\delta'_\Theta-\delta_\Theta)
(H_{w_0^\Theta\alpha})=
(-\lambda-\delta)(H_{w_0^\Theta\alpha})$.
Now it is clear that if $(\lambda+\delta)(H_\alpha)\neq 0$
for all $\alpha\in\triangle$ then $(\lambda'+\delta)(H_\alpha)\neq 0$
for all $\alpha\in\triangle$.
Since
$m=\#(\triangle^+\backslash<\Theta>),\
\lambda(H_\alpha)=\lambda'(H_\alpha)=0$ for $\alpha\in<\Theta>$ and
$\delta(H_\alpha)>0$ for $\alpha\in\triangle^+$, we get
$\#\{\alpha\in\triangle^+|(\lambda'+\delta)(H_\alpha)<0\}=
\#\{\alpha\in\triangle^+\backslash<\Theta>|(\lambda'
+\delta)(H_\alpha)<0\}=
\#\{\alpha\in\triangle^+\backslash<\Theta>|
(-\lambda-\delta)(H_{w_0^\Theta\alpha})<0\}=
\#\{\alpha\in\triangle^+\backslash<\Theta>|
(\lambda+\delta)(H_\alpha)>0\}=m-k$.
The Lemma is proved.

Retain the notations of Lemma 10 and assume that
$(\lambda+\delta)(H_\alpha)\neq 0$ for all $\alpha\in\triangle$.
It follows from Theorem 6 and Lemma 10 that the dual representations
$\tau$ and $\tau'$ of the algebra $g_c$ with
highest weights $\zeta$ and $\zeta'$ are realized in the cohomology
spaces $H^k(\M,{\cal S}L_\lambda)$ and $H^{m-k}(\M,{\cal S}L_{\lambda'})$
respectively. The spaces $H^k(\M,{\cal S}L_\lambda)$ and
$H^{m-k}(\M,{\cal S}L_{\lambda'})$ are dual. This is,
in fact, the Kodaira--Serre duality.

According to Theorem 7 and Lemma 9, the
$s$-modules ${\bf s}_\lambda$ and ${\bf s}_{\lambda'}$ are dual and the
associated function algebras $\A$ and $\A'$ on
$\M$ have representations $\rho$ and $\rho'$ in $H^k(\M,{\cal
S}L_\lambda)$ and $H^{m-k}(\M,{\cal S}L_{\lambda'})$, such that
$\tau=\rho\circ\sigma$ and $\tau'=\rho'\circ\sigma'$ respectively (here
all the notations have their usual meaning).
Let $n={\rm dim}\ \tau$.

{\bf Proposition 9. } {\it  For $f\in\A$ and $g\in\A'$ holds the equality
$$
	     n\int_\M fg\ d\mu={\rm tr}(\rho(f)(\rho'(g))^t).
$$

Proof. } Choose $u,v\in\U$ such that $f=\sigma_u,\ g=\sigma'_v$. Then,
using Eq. (7), one gets
$n\int fg\ d\mu=n\int \sigma_u\sigma'_v\ d\mu=
n\int \sigma_ul'_v1\ d\mu=n\int (l'_v)^t\sigma_u\ d\mu=
n\int l_{\check v}\sigma_u\ d\mu=n\int \sigma_{\check vu}\ d\mu=
{\rm tr}\ \tau(\check vu)={\rm tr}((\tau'(v))^t\tau(u))=
{\rm tr}(\rho(f)(\rho'(g))^t)$. The Proposition is proved.

{\bf 10. Covariant and contravariant symbols on flag manifolds}

Assume that $\lambda\in{\cal W}^\Theta$ is such that the (finite
dimensional) space ${\cal H}=H^0(\M,{\cal S}L_\lambda)$ of global
holomorphic sections of $L_\lambda$ is nontrivial. According to Theorem
6, this is the case iff $(\lambda+\delta)(H_\alpha)>0$ for all
$\alpha\in\triangle^+$ or, equivalently, iff $\lambda$ is a dominant
weight.

For any elements $q,q'$ of the same fiber of $L_\lambda^*$
denote by $h(q,q')$ their $K$-invariant
hermitian scalar product such that $h(q,q)=h(q)$.
Let $L^2(\M,L_\lambda)$ denote the Hilbert space of sections of
$L_\lambda$, square integrable with respect to the
$K$-invariant Hilbert norm $||\cdot||$ given by the formula
$||s||^2=\int_\M h(s)\ d\mu,\ s$ a section of $L_\lambda$.
Denote the corresponding hermitian scalar product
in $L^2(\M,L_\lambda)$ by $<\cdot,\cdot>$.

We introduce coherent states in ${\cal H}$ in a geometrically invariant
fashion, following \cite{C1}.  For $q\in L_\lambda^*$ the corresponding
coherent state $e_q$ is a unique element in ${\cal H}$ such that
the relation $<s,e_q>q=s\circ\pi(q)$ holds for all $s\in{\cal H}$. It is
known that the coherent states $e_q$ exist for all $q\in L^*_\lambda$ and
the mapping $L_\lambda^*\ni q\mapsto e_q\in{\cal H}$ is antiholomorphic.
For $c\in\C$ holds $e_{cq}=\bar c^{-1}e_q$.

The group $K$ acts on the sections of the line bundle $L_\lambda$ as
follows, $(ks)(x)=k(s(k^{-1}x))$ for $k\in K,\ x\in\M$ and $s$ a section
of $L_\lambda$. This action is unitary with respect to the scalar product
$<\cdot,\cdot>$. For any holomorphic section $s$ of $L_\lambda$ we have
$<ks,e_{kq}>kq=(ks)(kx)=k(s(x))$, therefore $<ks,e_{kq}>q=s(x)$.
On the other hand, $<ks,ke_q>=<s,e_q>=s(x)$, whence $ke_q=e_{kq}$.
The function $||e_q||^2h(q)$ is homogeneous of order 0 with respect
to $\C^*$-action and $K$-invariant. Thus it is identically constant. Set
$||e_q||^2h(q)=C$.

Let $A$ be an operator on ${\cal H}$. It is easy to check that the
function $\tilde f(q)=\\
<Ae_q,e_q>/<e_q,e_q>$ on the bundle $L_\lambda^*$
is constant on the fibers. Therefore there
exists a function $f_A$ on $\M$ such that $f_A\circ\pi=\tilde f$.

{\it Definition. }   Berezin's covariant symbol of an operator $A$
on ${\cal H}$ is the function $f_A$ on $\M$ given by the formula
$f_A(x)=<Ae_q,e_q>/<e_q,e_q>$ for any $q\in L^*_\lambda$ such that
$\pi(q)=x\in\M$.

The operator---symbol mapping $A\mapsto f_A$ is injective and thus
induces an algebra structure on the set of all
covariant symbols. The algebra of covariant symbols is isomorphic to
$End({\cal H})$.

Let $A$ be a holomorphic differential operator on $L_\lambda$.
Fix a local holomorphic trivialization $(U,s_0)$ of $L_\lambda$
and let $A_0$ denote the local expression of the operator $A$ on $U$.
Then for $x,y\in U$ we have
$<Ae_{s_0(y)},e_{s_0(x)}>s_0(x)=Ae_{s_0(y)}(x)=
s_0(x)A_0(e_{s_0(y)}(x)/s_0(x))$. The function $e_{s_0(y)}(x)/s_0(x)$
on $U\times U$ is holomorphic in $x$ and antiholomorphic in $y$.
Set $S(x)=<e_{s_0(x)},e_{s_0(x)}>=e_{s_0(x)}(x)/s_0(x)$.
Let $f$ be the covariant symbol of the operator $A$.
Since $A_0$ is a holomorphic differential operator on $U$, we get for
$q=s_0(x)$ that
$f(x)=<Ae_q,e_q>/<e_q,e_q>=(<Ae_{s_0(y)},e_{s_0(x)}>|_{y=x})/S(x)=
(A_0(e_{s_0(y)}(x)/s_0(x)))|_{y=x}/S(x)=A_0S(x)/S(x)$.

Introduce the function $\Phi=-\log h\circ s_0$ on $U$.
We have $S(x)=||e_{s_0(x)}||^2=C\exp\Phi$, whence $f(x)=A_0S(x)/S(x)=
e^{-\Phi}(A_0\ e^\Phi)=\check A1$, where $\check A$ is the
pushforward of the operator $A$ to $\M$. The formula $f=\check A1$
holds globally on $\M$.

For $u\in\U$ the pushforward to $\M$ of the operator $A_u$ on
$L_\lambda$ coincides with the operator $l_u$ related to the $s$-module
${\bf s}_\lambda$, $\check A_u=l_u$. Therefore the covariant symbol $f_u$ of
the operator $A_u$ on ${\cal H}$ can be expressed by the formula
$f_u=\check A_u1=l_u1=\sigma_u$. We have proved the following theorem.

{\bf Theorem 8. } {\it  Let $\lambda\in{\cal W}^\Theta$ be dominant.
Then the space ${\cal H}=H^0(\M,{\cal S}L_\lambda)$ of global holomorphic
sections of $L_\lambda$ is nontrivial. Endow it with the Hilbert space
structure via the norm $||\cdot||$ such that
$||s||^2=\int_\M h(s)\ d\mu,\ s\in{\cal H}$. Then for $u\in\U$ the
covariant symbol of the operator $A_u$ on ${\cal H}$ equals
$\sigma_u$, where $\sigma:\U\to C^\infty(\M)$ is the mapping associated
to the $s$-module ${\bf s}_\lambda$.}

According to Theorem 7, the representation $\rho$ of the algebra
$\A=\sigma(\U)$ in ${\cal H}$ is irreducible (the representation
$\tau=\rho\circ\sigma$ of the Lie algebra $g_c$ is irreducible with
highest weight $\lambda$).  Therefore any operator on ${\cal H}$ can be
represented as $A_u$ for some $u\in\U$.  Thus we get the following
corollary.

{\bf Corollary. } {\it The algebra $\A$ associated to the $s$-module
${\bf s}_\lambda$ coincides with the algebra of Berezin's covariant symbols of
the operators on ${\cal H}$.}

Let $P:L^2(\M,L_\lambda)\to{\cal H}$ be the orthogonal projection operator.
For a measurable function $f$ on $\M$ let $M_f$ denote the multiplication
operator by $f$. Introduce the operator $\hat f=PM_fP$ on ${\cal H}$.

{\it Definition. } A measurable function $f$ on $\M$ is called a
contravariant symbol of an operator $A$ on ${\cal H}$ if $A=\hat f$.

Let $s_1,s_2$ be holomorphic sections of $L_\lambda$.
Calculate the covariant symbol of the rank one operator
$A_0=s_1\otimes s_2^*$ in ${\cal H}$,
$$
f_{A_0}=\frac{<A_0e_q,e_q>}{<e_q,e_q>}
=\frac{<s_1,e_q><e_q,s_2>}{<e_q,e_q>}=
\frac{(s_1/q)\overline{(s_2/q)}}{||e_q||^2}=
\frac{h(s_1,s_2)}{||e_q||^2h(q)}.
$$
Since $||e_q||^2h(q)=C$, we obtain $f_{A_0}=h(s_1,s_2)/C$.
For any measurable function $g$ on $\M$
${\rm tr}(A_0\hat g)=<\hat
gs_1,s_2>=<gs_1,s_2>=\int h(gs_1,s_2)\ d\mu= \int gh(s_1,s_2)\
d\mu=C\int f_{A_0}g\ d\mu$. Therefore for any operator $A$ on ${\cal H}$
holds ${\rm tr}(A\hat g)=C\int f_Ag\ d\mu$. Taking $A={\bf 1},\ g=1$
we immediately obtain that $C=n={\rm dim}\ {\cal H}$.

{\bf Proposition 10.} {\it A measurable function $g$ on $\M$ is a
contravariant symbol of an operator $B$ on ${\cal H}$ iff for any
operator $A$ on ${\cal H}$ holds the formula
${\rm tr}(AB)=n\int f_Ag\ d\mu$.}

The proof is straightforward.

Let $\lambda\in{\cal W}^\Theta$ be
dominant. Set $\lambda'=-\lambda-2\delta'_\Theta$. For the $s$-module
${\bf s}_{\lambda'}$ dual to ${\bf s}_\lambda$ let $\sigma':\U\to C^\infty(\M)$
denote the mapping associated to ${\bf s}_{\lambda'}$, $\rho'$ be the
corresponding representation of the algebra ${\cal A'}=\sigma'(\U)$ in
$H^m(\M,{\cal S}L_{\lambda'})$ The spaces
${\cal H}=H^0(\M,{\cal S}L_\lambda)$ and
$H^m(\M,{\cal S}L_{\lambda'})$ are dual as representation
spaces of the group $G$.
The following theorem is a direct consequence of Propositions  9,10
and Theorem 8.

{\bf Theorem 9. } {\it A function $f\in {\cal A'}$ is a contravariant
symbol of the operator $(\rho'(f))^t$ in ${\cal H}$.}

{\bf 11. Quantization on flag manifolds}

Now we are ready to put together various results obtained above to
give examples of quantization on a generalized flag manifold $\M$
endowed with a pseudo-K\"ahler metrics.

Let $\lambda\in{\cal W}^\Theta$ be such that for all
$\alpha\in\triangle\backslash<\Theta>$ holds $\lambda(H_\alpha)\neq 0$.
According to Proposition 8, in this case  the
$s$-module ${\bf s}_\lambda$ on $\M$
is nondegenerate. Denote by $\omega$ the 2-form associated to
${\bf s}_\lambda$. This form is pseudo-K\"ahler, and the index of inertia
of the corresponding pseudo-K\"ahler metrics equals
$l=\#\{\alpha\in\triangle^+\backslash<\Theta>|\lambda(H_\alpha)<0\}$.
Denote by $\A_\hbar$ the algebra of functions on $\M$ associated to the
$s$-module $(1/\hbar){\bf s}_\lambda$. It follows from Theorem 5
that any functions $f,g\in F(\M)$ belong to $\A_\hbar$ for all but a
finite number of values of $\hbar$. The product $f*_\hbar g$ expands to
the uniformly and absolutely convergent Taylor series in $\hbar$ at the
point $\hbar=0$, $f*_\hbar g=\sum_{r\geq 0} \hbar^r C_r(f,g)$, where
$C_r(\cdot,\cdot), r=0,1,\dots$, are bidifferential operators which
define the deformation quantization with separation of variables on $\M$
corresponding to the (non-deformed) pseudo-K\"ahler form $\omega$.

For $n\in{\cal N}$ holds $n\lambda\in{\cal W}^\Theta$.
Theorem 7 implies that for $\hbar=1/n$ the algebra $\A_\hbar$
has a natural geometric representation $\rho_\hbar$ in the
sheaf cohomology space of the line bundle $L_{n\lambda}=(L_\lambda)^n$,
${\cal H}_\hbar=H^{k_n}(\M,{\cal S}L_{n\lambda})$, where
$k_n=\#\{\alpha\in\triangle^+|(n\lambda+\delta)(H_\alpha)<0\}$.
Since for $\alpha\in\triangle^+$ holds $\delta(H_\alpha)>0$,
only those $\alpha\in\triangle^+$ contribute to $k_n$ for which
$\lambda(H_\alpha)<0$. Therefore $k_n=l$ for $n>>0$. In other words,
for sufficiently small values of $\hbar=1/n$ the dimension of the
sheaf cohomology the representation $\rho_\hbar$ is realized in is equal
to the index of inertia $l$ of the pseudo-K\"ahler metrics on
$\M$ corresponding to the (1,1)-form $\omega$.

We have obtained pseudo-K\"ahler quantization on a generalized flag
manifold.

Now assume $\lambda\in{\cal W}^\Theta$ is dominant in the rest of the
paper. Then the metrics corresponding to the (1,1)-form $\omega$ on $\M$
is positive definite, i.e., $\omega$  is a K\"ahler form, and Theorem 8
implies that for $\hbar=1/n$ the space
${\cal H}_\hbar=H^0(\M,{\cal S}L_{n\lambda})$ is the space of
global holomorphic sections of the line bundle
$L_{n\lambda}=(L_\lambda)^n$ and
$\A_\hbar$ is the corresponding algebra of Berezin's covariant symbols on
$\M$.  Thus we arrive at Berezin's K\"ahler quantization on $\M$ and
identify the associated formal deformation quantization obtained in
\cite{M}, \cite{C2} with the quantization with separation of variables,
corresponding to the non-deformed K\"ahler form $\omega$.

Consider the $s$-module ${\bf s}(\hbar)=
-{\bf s}_\lambda+\hbar{\bf s}_{can}$. It depends rationally on $\hbar$
and is regular at $\hbar=0$. Denote by $\omega_{can}$ the (1,1)-form
associated to the canonical $s$-module ${\bf s}_{can}$ on $\M$. Then
the form, associated to ${\bf s}(\hbar)$ is $-\omega+\hbar\omega_{can}$.
The $s$-module $(1/\hbar){\bf s}(\hbar)$ is dual to $(1/\hbar){\bf
s}_\lambda$. Denote by $(\A'_\hbar,*'_\hbar)$ the algebra of functions on
$\M$ associated to the $s$-module $(1/\hbar){\bf s}(\hbar)$.  For any
functions $f,g\in F(\M)$ the product $f*'_\hbar g$ depends rationally on
$\hbar$ and is regular at $\hbar=0$. The asymptotic expansion
of the product $f*'_\hbar g$ gives rise to the deformation
quantization with separation of variables $({\cal F},\star')$
corresponding to the formal deformation of the negative-definite K\"ahler
form $-\omega$, $\omega'=-\omega+\nu\omega_{can}$.

If $\hbar=1/n$ then the algebra $\A'_\hbar$ has a representation
$\rho'_\hbar$ in the space $H^m(\M,{\cal S}L_{\lambda'_n})$,
where $m={\rm dim}_\C\M$ and $\lambda'_n=-n\lambda-2\delta'_\Theta$.
The space $H^m(\M,{\cal S}L_{\lambda'_n})$ is dual to ${\cal H}_\hbar$
and Theorem 9 implies that any function $f\in\A'_\hbar$ is a
contravariant symbol of the operator $(\rho'_\hbar(f))^t$ in
the space ${\cal H}_\hbar$. The mapping
$\A'_\hbar\ni f\mapsto (\rho'_\hbar(f))^t$ is an anti-homomorphism.
Thus, in order to obtain quantization on $\M$ by contravariant symbols
(it is usually called Berezin-Toeplitz quantization, see \cite{Sch}),
we have to consider the algebras $(\tilde\A_\hbar,\tilde *_\hbar)$,
opposite to $(\A'_\hbar,*'_\hbar)$. Then
$\tilde\A_\hbar\ni f\mapsto (\rho'_\hbar(f))^t$ will be a representation
of the algebra $\tilde\A_\hbar$. The corresponding deformation
quantization $({\cal F},\tilde\star)$ is opposite to $({\cal F},\star')$.
As it was shown in Section 5, this quantization is also a
quantization with separation of variables, though with respect to the
opposite complex structure on $\M$. It corresponds to the formal
(1,1)-form  $-\omega'=\omega-\nu\omega_{can}$ on the opposite complex
manifold $\bar\M$. This form is a formal deformation of the K\"ahler form
$\omega$ on $\bar\M$. (The metrics on $\bar\M$, corresponding to the
(1,1)-form $\omega$ is a negative-definite K\"ahler metrics.)

It would be interesting to compare the deformation quantization
associated to Berezin-Toeplitz quantization on a general compact
K\"ahler manifold in \cite{Sch} with deformation quantization with
separation of variables.

{\it Acknowledgements.} I wish to express my deep gratitude to
Professor M.S. Narasimhan for inviting me to the
International Centre for Theoretical Physics (ICTP)
and to the ICTP for their warm hospitality.

\end{document}